\documentclass{amsart}

\usepackage{amsmath,amssymb,amsthm, mathrsfs, mathtools, bm, amsfonts}
\usepackage{mathabx}
\usepackage{amscd,mathtools}
\usepackage[utf8]{inputenc}
\usepackage[english]{babel}
\usepackage{cite}
\usepackage{color}
\usepackage[pagebackref,colorlinks,citecolor=blue,linkcolor=red]{hyperref}

\usepackage{float}
\usepackage{tikz}
\usepackage{lmodern}
\usetikzlibrary{decorations.pathmorphing}

\newcommand{\E}{{\mathcal E}}

\newcommand{\PSL}{{\rm{PSL}}}

\usepackage{multicol}
\tikzset{snake it/.style={decorate, decoration=snake}}
\usetikzlibrary{automata}
\usetikzlibrary{cd}

\usepackage{comment}
\usepackage{mathtools}

\usepackage{csquotes}

\usepackage{amsthm}

\usepackage{cleveref}

\usepackage{longtable}

\newtheorem{theorem}{Theorem}[section]
  \newtheorem{proposition}[theorem]{Proposition}

  \newtheorem{corollary}[theorem]{Corollary}
  \newtheorem{lemma}[theorem]{Lemma}
  \newtheorem{claim}[theorem]{Claim}

\theoremstyle{definition}
\newtheorem{defn}{Definition}

\newtheorem{remark}[theorem]{Remark}

\newcommand{\rank}{{\rm{rank}}}
\newcommand{\sys}{{\rm{sys}}}
\newcommand{\esssys}{{\rm{sys_e^*}}}
\newcommand{\combsys}{{\rm combsys}}
\newcommand{\Tet}{{\rm {\frak Tet}}}
\newcommand{\T}{{\mathcal{T}}}
\newcommand{\F}{{\mathcal{F}}}
\newcommand{\CAT}{{\rm CAT}}
\newcommand{\sss}{\hspace{-.05cm} \smallsetminus \hspace{-.05cm}}

\usepackage{xcolor}
\usepackage{fancyvrb}
\definecolor{Green}{RGB}{0,180,60}
\definecolor{Purple}{RGB}{150,0,250}
\definecolor{mypink}{RGB}{255, 0, 255}

\newcommand{\yandi}[1]{\marginpar{\tiny \color{blue} Yandi: #1}}

\newcommand{\blue}[1]{\textcolor{blue}{#1}}

\usepackage{listings}

\usepackage{xparse}

\NewDocumentCommand{\gap}{v}{%
\textbf{\texttt{\textcolor{blue}{#1}}}%
}
\NewDocumentCommand{\codeinput}{v}{%
\texttt{\textcolor{red}{#1}%
}}
\NewDocumentCommand{\codeoutput}{v}{\texttt{\textcolor{black}{#1}}}
\NewDocumentCommand{\codecomment}{v}{\texttt{\textcolor{purple}{#1}}}
\usepackage{array} 

\DeclareMathOperator{\aux}{aux}

\title{Filling knots and essential geodesic systole}

\author{Christopher J Leininger}\thanks{CJL was supported by NSF grant DMS-2305286}
\address{University of Notre Dame}
\email{cleining@nd.edu}
\urladdr{\url{https://sites.google.com/view/chris-leiningers-webpage/home}}

\author{Yandi Wu}
\address{Rice University}
\email{yw220@rice.edu}
\urladdr{\url{https://yandiwu.github.io/}}

\begin{document}

\maketitle 

\begin{abstract} We prove that every closed, orientable $3$--manifold contains hyperbolic knots with arbitrarily large essential geodesic systole. That is, for any $3$--manifold $M$, and $R> 0$ there is a hyperbolic knot $K \subset M$ such that every closed geodesic in the knot complement which is non-null-homotopic in $M$ has length at least $R$.  The construction also produces filling knots, in the sense of Freedman-Krushkal \cite{FreeKrus2023}.  This answers two questions of Biringer.
\end{abstract}

\section{Introduction}
The \emph{systole} of a Riemannian manifold $N$, $\sys(N)$, is the length of the shortest closed geodesic in $N$. 
In hyperbolic geometry, there are interesting relationships between the systole and the topology of the manifold. For example, if $N$ is a closed \emph{hyperbolic} 3-manifold, White found that $\sys(M)$ is bounded by a function of the \emph{rank}; i.e.~the minimal number of generators for $\pi_1(M)$ \cite{White}; see also Kapovich-Weidmann \cite{kapovichweidmann}. 

Another interesting connection (relevant for the current paper) was proved by Adams and Reid in \cite{adamsreid} for link complements.  They showed that if $M$ is {\em any} closed $3$--manifold which does not admit a hyperbolic metric, then for {\em any} link $L \subset M$ such that $M \sss L$ is hyperbolic, one has $\sys(M \sss L)\leq 7.35534...$.  The hypothesis that $M$ admits no hyperbolic metric is necessary.  On the other hand, for a closed hyperbolic $3$--manifold, Lakeland and the first author found that for {\em any} link $L \subset M$ with $M \sss L$ hyperbolic, $\sys(M \sss L)$ is bounded above by a function of the volume of $M$ \cite{lakelandleininger}. In either case, for a fixed $M$, there a uniform upper bound on $\sys(M \sss L)$ that does not depend on the link $L$ itself. 

We note one similarity between the proofs of \cite{adamsreid} and \cite{lakelandleininger}: the short geodesics which illustrate the bound on $\sys(M \sss L)$ are null-homotopic in the ambient manifold, $M$.
This inspires the following definition. 

\begin{defn}[Essential geodesic systole]\label{def:esssys}
The \emph{essential geodesic systole} of a link $L \subset M$ is the length of the shortest closed geodesic in $M \sss L$ which is not null-homotopic in $M$,
\[ \esssys(K) = \inf\{\ell(\gamma) \mid \gamma \mbox{ a closed geodesic in } M \sss L \mbox{ with } i \circ \gamma \mbox{ non-null-homotopic in } M \}. \]
\end{defn}

In an earlier paper, the authors of this paper proved the following: 

\begin{theorem}[Theorem~1.2,\cite{lw}]
\label{thm:esssys}
Let $M$ be a closed, orientable $3$-manifold with $\pi_1(M) \neq \{\ast\}$. Then for every $R > 0$, there exists some hyperbolic link $L \subset M$ with \[\esssys(L) > R.\]
\end{theorem}

\begin{remark}\label{rem: esssys new defn} \Cref{def:esssys} differs from the definition of \emph{essential systole}, denoted $\sys_e(L)$, given in \cite{lw}; there, the infimum was taken over all closed {\em curves} rather than closed geodesics, although \Cref{thm:esssys} holds for both definitions (see \cite{lw} Remark 1.3).  We use this new definition since it allows the knot itself to be non-null-homotopic in $M$, which will be more convenient for the current paper. 
Thus, when we say \emph{essential systole} in the present paper, we implicitly mean essential geodesic systole.   
\end{remark}

The purpose of this paper is to strengthen the theorem above, as follows, answering a question of Biringer for essential geodesic systole (see \cite[Question 4.1]{lw}).

\begin{theorem}\label{thm:longessentialsys} Let $M$ be an orientable $3$--manifold. Let $R > 0$. Then there exists some hyperbolic knot $K \subset M$ such that 
\[\esssys(K) > R.\]
\end{theorem} 

In \cite{lw}, we used \Cref{thm:esssys} to prove the existence of \emph{filling links} in any closed, orientable 3-manifold with nontrivial fundamental group. A filling link, originally defined by Freedman and Krushkal in \cite{FreeKrus2023}, is a link that ``cuts through all the topology" of its ambient 3-manifold by killing the relations in \textit{any} group presentation of $\pi_1(M)$ with respect to any minimal generating set, {\em viewed inside $M \sss L$}. 

More precisely, for a closed $3$--manifold, $M$, a {\em $1$--spine} (also called a {\em carrier graph}) for $M$ is a graph $\Gamma$ with $\rank(\pi_1(\Gamma)) = \rank(\pi_1(M))$ together with a map $f \colon \Gamma \to M$ such that $f_* \colon \pi_1(\Gamma) \to \pi_1(M)$ is surjective.
Given a link $L \subset M$, let $i \colon M \sss L \to M$ be the inclusion.  We say that $f \colon \Gamma \to M \sss L$ is an $L$--relative $1$--spine if $i \circ f \colon \Gamma \to M$ is a $1$--spine.
Freedman and Krushkal declare that a link $L \subset M$ is a {\em filling link} if for every $L$--relative $1$--spine $f \colon \Gamma \to M \sss L$, we have
\[ f_*\colon \pi_1(\Gamma) \to \pi_1(M \sss L) \]
is injective. They then asked in \cite{FreeKrus2023} whether filling links exist in every closed $3$--manifold, and in \cite{lw}, we proved that the answer is `yes' for all closed, orientable $3$--manifolds. 

After \cite{lw}, Biringer asked the authors whether filling {\em knots} exist (see \cite[Question 4.1]{lw}), and we again prove that the answer is `yes'.




\begin{theorem} \label{T:filling knot}
    Any closed, orientable $3$--manifold $M$ contains a filling knot \\$K \subset M$.
\end{theorem}

We recall that the filling links with large essential systole from \cite{lw} are constructed using sequences of triangulations of the ambient manifold $M$, and the number of link components grows with the number of tetrahedra in the associated triangulation.  Thus, the construction of knots with large essential geodesic systole is significantly harder and requires new ideas.  Moreover, (as noted above) in \cite{lw}, we prove that the links are filling as a direct consequence of having large essential systole, applying the main result of \cite{kapovichweidmann}.  The distinction between large essential systole and large essential {\em geodesic} systole (see \Cref{rem: esssys new defn}), means that the existence of filling knots in the current paper also requires new ideas.

\begin{remark} As in \cite{lw} (also as suggested by Biringer), for any $m \in \mathbb{N}$, we prove the existence of a \emph{full rank-$m$ filling} knot, which is a stronger version of \Cref{T:filling knot}. When $m = \pi_1(M)$, a full rank-$m$ filling knot is also filling. However, for ease of exposition, we defer this technicality to \Cref{sec:constructingknots}; see \Cref{thm:main2specific}. Similarly, there is a more precise version of \Cref{thm:longessentialsys}; see \Cref{thm:main1specific}.  
\end{remark}

We now give a brief outline of the paper. In \Cref{sec:prelims}, we include some important tools used in the rest of the paper. In \Cref{section:l_prime_construction}, given a manifold $M$, we construct a family of knots used to prove both main theorems. First, we build ``augmented" versions of the links seen in \cite{lw} by replacing certain totally geodesic subsurfaces of the original link complements with hyperbolic tangle complements, which we call \emph{connectors}. These newly glued-in connectors allow us ``join" together link components by applying half-twists to their totally geodesic subsurfaces; c.f.~\cite{Adams-augment,Lackenby-augment,Purcell-augment}. We then apply large Dehn fillings and surgeries to reduce the number of link components to just $1$. 

Finally, in Sections \ref{sec:longessentialsystole} and \ref{sec:fillingknots}, we prove Theorems \ref{thm:longessentialsys} and \ref{T:filling knot} respectively. The proofs use similar approaches, but slightly different tools. The existence of knots with large essential geodesic systole relies on Thurston's Hyperbolic Dehn Surgery Theorem and other tools from hyperbolic geometry. While it is probably true that these hyperbolic knots are also filling, it is more natural and convenient to endow our knot complements with negatively-curved metrics. Thus, we instead prove the existence of filling links using locally CAT($-1/2$) metrics coming from negatively-curved Dehn fillings derived from the Gromov--Thurston $2\pi$ Theorem.\\  

\noindent \textbf{Acknowledgments.} The authors are very grateful to Alan Reid for many useful conversations about the project, and especially for suggesting the proof of \Cref{thm:connectorhyperbolic} using the \lstinline{verify_hyperbolicity} function in SnapPy. They also thank Ian Biringer for initially suggesting the problem, and for helpful discussions during the second author's visit to Boston College.\\

\noindent \textbf{AI disclosure.}
No AI was used for any part of this work, other than that imposed by {\em Google} in its search engine.

\section{Preliminaries}
\label{sec:prelims}

We recall some basic tools used in our constructions in the rest of this paper.
We let $\mathbb H^3$ denote hyperbolic $3$--space with metric $d_{\mathbb H^3}$.  In the upper-half space model, the orientation preserving isometry group is $\PSL_2(\mathbb C)$ acting by extension of M\"obius transformation.
A(n orientable) \textit{hyperbolic 3-manifold} is the quotient of $\mathbb{H}^3$ by a discrete, torsion-free subgroup of $\PSL_2(\mathbb{C})$. Given a closed $3$-manifold $M$, we say that a link $L \subset M$ is \emph{hyperbolic} if $M \sss L$ admits a complete hyperbolic structure of finite volume. See \cite{thurston82,Thurstonnotes} for more on hyperbolic geometry and hyperbolic manifolds.



Recall that an $n$-manifold $X$ with a geodesic metric is \emph{locally $\CAT(-1)$} if every point has a neighborhood in which geodesic triangles are at least as thin as those in the hyperbolic plane, $\mathbb H^2$; if the metric is complete, then by the Cartan-Hadamard Theorem \cite[Theorem~II.4.1]{bh}, its universal cover $\tilde{X}$ with its pull-back metric is $\CAT(-1)$ ({\em all} triangles satisfy the above comparison).  For $\kappa < 0$, the (locally) $\CAT(\kappa)$ condition is similarly defined, by comparing with triangles in $\mathbb H^2$ with the metric scaled by $\sqrt{|\kappa|}$.   In this paper, we will make abundant use of the following consequence of Thurston's Hyperbolization Theorem \cite{thurston82,Thurston-hyp,Morgan,McMullen}, and which allows us to endow the locally $\CAT(-1)$ link complements constructed in this paper with hyperbolic metrics.  These locally $\CAT(-1)$ and hyperbolic link complements will have the added property that their universal covers are quasi-isometric.  A proof, appealing to Hyperbolization, is sketched in \cite{lw}.

\begin{theorem}[\!\cite{lw} Theorem 2.2] \label{T:CAT-1 to hyp}
Suppose a 3--manifold $N$ admits a complete, finite-volume, locally $\CAT(-1)$, piecewise hyperbolic metric which is hyperbolic in the cusps. Then $N$ admits a complete, finite-volume hyperbolic metric.  Moreover, this metric can be chosen so that the identity on the universal covers is a quasi-isometry.
\end{theorem}

\subsection{Thick-thin decomposition} The \emph{$\epsilon$-thin (resp. thick) part} of a manifold $M$, denoted $M_{(0,\epsilon)}$ (resp. $M_{[\epsilon, \infty)}$), is the set of points $x \in M$ such that the injectivity radius of $M$ at $x$ is less than (resp. greater than or equal to) $\epsilon$. One important consequence of the Margulis Lemma is that there is a value $\epsilon_0 > 0$, known as the ($3$-dimensional) \emph{Margulis constant}, such that for any $\epsilon \leq \epsilon_0$ and $3$-manifold $M$, $M_{(0, \epsilon)}$ is either a \emph{Margulis tube}, the tubular neighborhood of a short curve, or a cusp neighborhood; see e.g.~\cite[Chapter~5]{Thurstonnotes}.  

Another consequence of the Margulis Lemma is the following, which is an important ingredient in our arguments.
Given a discrete, torsion--free subgroup $\Gamma < \PSL(2,\mathbb C)$, the $\epsilon$--thin part of $\mathbb H^3$ is the preimage of the $\epsilon$--thin part of $M/\Gamma$.  For $\epsilon > \epsilon_0$, the $\epsilon$--thin part is a disjoint union of horoballs and neighborhoods of geodesics.  Note that the $\epsilon$--thin part is precisely the set of points which are moved a distance less than $2\epsilon$ by some nontrivial element of $\Gamma$.


\begin{lemma} \label{lem:bounded abelian}
 Given $K > 0$ and $0 < \epsilon < \epsilon_0$, there exists $\eta_0 > 0$ with the following property.  Suppose $\Gamma < \PSL(2,\mathbb C)$ is a discrete group and $\widetilde z_1,\ldots,\widetilde z_\eta \in \mathbb H^3$ are points in distinct $\epsilon$--thin parts for $\Gamma$.  If there is some $\widetilde z_0 \in \mathbb H^3$ so that $d_{\mathbb H^3}(\widetilde z_0,\widetilde z_i) \leq K$ for all $i\geq 1$, then $\eta < \eta_0$.
\end{lemma}
\begin{proof}
    Since $0 < \epsilon < \epsilon_0$, then for $\delta = \tfrac{\epsilon_0 - \epsilon}2$, the ball $B_i$ of radius $\delta$ about $\widetilde z_i$ is contained in the $\epsilon_0$--thin part for $\Gamma$.  Moreover, $\{B_i\}_i$ is a collection of pairwise disjoint balls of radius $\delta$ contained in the ball of radius $K+\delta$.  Thus we can take $\eta_0$ to be an integer larger than the ratio of the volume of a ball of radius $K+\delta$ to that of a ball of radius $\delta$.
\end{proof}

Another key ingredient in this paper is the fact that geodesics intersecting thick parts of 3-manifolds that enter Margulis tubes or cusps are long. This fact is made precise by the following lemma. 

\begin{lemma}[\cite{bm, meyerhoff}] \label{lem:thinpart} There is a continuous positive function \[R: [0, \epsilon_0] \rightarrow \mathbb{R}_{\geq 0} \text{ such that } \lim\limits_{\epsilon \rightarrow 0} R(\epsilon) = \infty\] with the following property.  Given $\epsilon \leq \epsilon_0$, let $T_{\epsilon} \subset T_{\epsilon_0}$ be connected components of the $\epsilon$-- and $\epsilon_0$--thin parts, respectively, of a hyperbolic manifold $M$. Then
\[d_{M}\big(\partial T_{\epsilon}, \partial T_{\epsilon_0}\big) \geq R(\epsilon).\] 
\end{lemma}

\subsection{Dehn filling/surgery} A crucial ingredient in our construction of knots rather than links is Thurston's hyperbolic Dehn filling theorem; see \cite[Chapter 5]{Thurstonnotes} and \cite[Chapter E]{BenedettiPetronio}.  This allows cusps of a finite volume hyperbolic $3$--manifold to be ``filled in" by solid tori while preserving hyperbolicity and approximating some geometric features of the original manifold. 

Suppose $M$ is a compact $3$-manifold and $\{T_1, \ldots, T_k\} \subset \partial M$ its boundary components, each of which is homeomorphic to a torus. Let $\alpha_j \subset T_j$ be an essential simple closed curve in $T_j$. Let $U_j$ be a solid torus, with meridian $\beta_j$. Consider a homeomorphism \[\phi_j: (\partial U_j, \beta_j) \rightarrow (T_j, \alpha_j)\] which sends each $\beta_j$ homeomorphically to each $\alpha_j$. The resulting manifold obtained by ``filling in" each $T_j$ via the map $\phi_j$ (that is, gluing $U_j$ to $T_j$) is the \emph{Dehn filling of $M$ along $\alpha_j$}. Up to homeomorphism, the result of Dehn filling depends only on the isotopy class of $\alpha_j$ in $T_j$.  Such an isotopy class is called a {\em slope} on $T_j$, and we often use the same notation for the simple closed curve and its isotopy class. We can parameterize the isotopy classes $\alpha_j$ by considering framings of each $T_j$ by a choice of longitude $\mu_j$ and meridian $\nu_j$ (i.e.~a basis for $\pi_1(T_j)$). We can encode $\alpha_j$ by the tuple $(p_j, q_j) \in \mathbb{R}^2$ (up to sign) in $\pi_1(T_j)$, if (after orienting $\alpha_j$,) $\alpha_j = p_j\mu_j + q_j\nu_j$. Setting $u_j = (p_j, q_j)$, we write both $M(u_1,\ldots, u_k) = M(\alpha_1,\ldots,\alpha_k)$ as the manifold obtained from Dehn filling along each $\alpha_j$.  We also allow $\alpha_j$ or $u_j$ to be ``$\infty$", by which we mean that we perform no Dehn filling on the boundary component.


If $M$ is the interior of a compact manifold with torus boundary, we similarly define Dehn filling by deleting product neighborhoods of the ends to produce a compact $3$--manifold whose boundary components are all homeomorphic to tori; these tori, viewed in the original non-compact manifold $M$, will be called {\em peripheral tori}.  If a filling slope is $\infty$, then we glue the product neighborhood of the end back in.  If $L \subset M$ is a (tame) link in compact $3$--manifold $M$, then a {\em Dehn surgery} on $L$ is the process of drilling out $L$ (i.e.~removing $L$), and performing Dehn filling on the resulting open manifold (or equivalently, removing tubular neighborhoods of $L$ and Dehn filling torus boundary components of the resulting manifold).




Let $H$ be the space of all finite-volume 3--manifolds equipped with the \emph{geometric topology}.  That is, a sequence of 3--manifolds $\{N_i\} \subset H$ converges to $N \in H$ if and only if for some sequence of positive numbers $\epsilon_i \rightarrow 0$, there is a sequence $\{\phi_i\}$, where each $\phi_i$ is an $(1 + \epsilon_i)$-biLipschitz diffeomorphism between the $\epsilon_i$-thick part of $N_i$ and the $\epsilon_i$-thick part of $N$. We are now ready to state the following theorem due to Thurston; see \cite[Theorem 5.8.2]{Thurstonnotes}, \cite[Sections E.5-E.6]{BenedettiPetronio}.

\begin{theorem}[Hyperbolic Dehn Surgery Theorem \cite{Thurstonnotes}] \label{thm:dehnfilling} Suppose $N$ is a finite volume, cusped hyperbolic $3$--manifold.  Then $N(u_1,\ldots,u_k)$ is hyperbolic, provided each $u_i$ avoids a finite set of (non-$\infty$) slopes on the $i^{th}$ cusp.  Furthermore, if $u_{i,n} \to \infty$ in $\mathbb R^2 \cup \{\infty\}$ for each $i$, then $\{N(u_{1,n},\ldots,u_{k,n})\}_n$ converges to $N$.
\end{theorem}

\subsection{A theorem of Kapovich--Weidmann} Lastly, we introduce a result which will be used to show that the knots we construct are filling knots. Let $G$ be a free group of rank $n$. Two generating sets for $G$ are \emph{Nielsen equivalent} if they differ by a sequence of \emph{Nielsen transformations}, which are analogs of elementary row operations on $GL_n(\mathbb{Z})$ (see \cite{nielsen} for details), and which generate the automorphism group of $G$. In other words, Nielsen transformations map ordered an (ordered) basis of $G$ to another basis of $G$. Like in \cite{lw}, in \Cref{sec:fillingknots}, we will make use of the following version of a theorem due to Kapovich--Weidmann, which roughly states that a group that acts by isometries on a $\delta$-hyperbolic space $(\delta > 0)$ is either free or contains a nontrivial element with short translation length. 

\begin{theorem}[\!\cite{kapovichweidmann}] \label{thm:KW} For any integer $n > 0$, there exists a constant $C(n)$ with the following property. Suppose a group $G =  \langle g_1,\ldots,g_n\rangle$ acts by isometries on a $\delta$-hyperbolic geodesic metric space $(X, d)$, where $\delta > 0$.  Then one of the
following holds:
\begin{enumerate}
\item The group $G$ is free with basis $(g_1, g_2,\ldots, g_n)$, or 
\item The tuple $(g_1, g_2, \ldots, g_n)$ is Nielsen-equivalent to a tuple $(g_1', g_2', \ldots, g_n')$ such that for some $x \in X$, $d_X(g'_1.x, x) < \delta C(n)$.
\end{enumerate}
\end{theorem}
\section{Constructions of links}
\label{section:l_prime_construction}

We fix a closed, orientable $3$--manifold $M$ for the remainder of this paper.  We first recall the construction of links $\{L\}$ in $M$ with large essential systole, built from special triangulations of $M$, and remind the reader of a few relevant facts about them.  We then define a new tangle, and describe features of a nice hyperbolic metric on it.  Next, we explain how to modify $L$ to obtain a new link $L'$ in $M$ which will serve as the primary ingredient, via Dehn surgery, for the construction of knots that prove the main theorems.  We end with the construction of the knots and state more precise versions of the theorems from the introduction.

\subsection{Triangulations and combinatorial systole} \label{S:triangulations and combsys}

To describe the construction of the knots in a closed $3$--manifold $M$ needed to prove both Theorem~\ref{thm:longessentialsys} and Theorem~\ref{T:filling knot}, we recall two definitions from \cite{lw}.

Given a triangulation $\tau$ of $M$, we define the {\em combinatorial systole} of $\tau$ to be the length of the shortest non-null-homotopic loop in $M$, measured with respect to the graph metric on the $1$--skeleton of $\tau$.  There are various ways to make this precise, and the one that is most useful for our purposes is to consider the subgraph $\hat \tau$ of the $1$--skeleton of the barycentric subdivision of $\tau$ consisting of all $1$--simplices that meet the barycenters of the tetrahedra of $\tau$.  This is a bipartite graph, with vertices partitioned into barycenters of tetrahedra and all other barycenters.  A pair of tetrahedra-vertices are distance two apart if and only if the tetrahedra intersect.  The inclusion induces a surjective homomorphim $\pi_1(\hat \tau) \to \pi_1(M)$, and we define
\[ \combsys(\tau) = \min\{\ell_{\hat{\tau}}(\gamma) \mid \gamma \mbox{ a loop in }\hat \tau, \gamma \not\simeq \ast \mbox{ in }M \}. \]

We will utilize a special type of triangulation $\tau$ of our ambient manifold $M$ which we call a {\em Cooper-Thurston triangulation}.  These were defined by Cooper and Thurston in \cite{CoopThur1988}, and have the property that the link of every vertex is one of $5$ particular $2$--complexes (each of which is flag), and each edge has $4$, $6$, $8$, or $10$ tetrahedra glued around it.  In \cite[Proposition~3.2]{lw}, we proved the following.
\begin{proposition}
    For any closed, orientable $3$--manifold $M$ with $\pi_1M \neq \{\ast\}$ and $r > 0$, there exists a Cooper-Thurston triangulation $\tau$ with $\combsys(\tau) > r$. 
\end{proposition}

\subsection{Link construction and metrics} \label{subsection:fillinglinkrecap} 

Given a Cooper-Thurston triangulation $\tau$ of $M$, in \cite{lw} we constructed an associated hyperbolic link $L = L(\tau) \subset M$ whose complement $M_L = M \sss L$ admits an explicit, complete, locally $\CAT(-1)$, singular hyperbolic metric, with singular locus contained in the $1$--skeleton of $\tau$.  In particular, by \Cref{T:CAT-1 to hyp}, the link is itself hyperbolic, and we write $M_L^h$ for the hyperbolic metric coming from that theorem. The essential systole of $L$ is bounded below by a function of the combinatorial systole of $\tau$ (see \cite{lw}).  Before we describe the new constructions of links and knots, we first recall the construction from \cite{lw}.

\begin{figure}[h]
\begin{tikzpicture}
    \begin{scope} [scale = 2]
        \draw[thick] (0,0) -- (1,-1) -- (2,0) -- (1,1) -- (0,0) -- (1,-1);
        \draw[thick] (0,0) -- (.5,0);
        \draw[thick] (.75,0) -- (.95,0);
        \draw[thick] (1.05,0) -- (1.25,0);
        \draw[thick] (1.5,0) -- (2,0);
        \draw[thick] (1,-1) -- (1,1);
        \draw[ultra thick, blue] (1,.2) .. controls (.7,.1) and (.5,0) .. (.4,-.4);
        \draw[ultra thick, violet] (1,-.2) .. controls (1.3,-.1) and (1.5,0) .. (1.6,.4);
        \node at (.5,-.7) {$\frac{\pi}2$};
        \node at (1.1,.5) {$\frac{\pi}2$};
        \node at (1.5,.7) {$\frac{\pi}2$};
        \node at (1.5,-.7) {$\frac{\pi}3$};
        \node at (.5,.7) {$\frac{\pi}4$};
        \node at (1.3,-.3) {$\frac{\pi}3$};
        \draw[->,>=stealth] (1.4,-.3) .. controls (1.5,-.3) and (1.7,-.2) .. (1.7,-.05);
    \end{scope}
    \begin{scope}[shift = {(5,-.5)},scale=1.5]
    \node at (2,.4) {\includegraphics[width=0.3\linewidth]{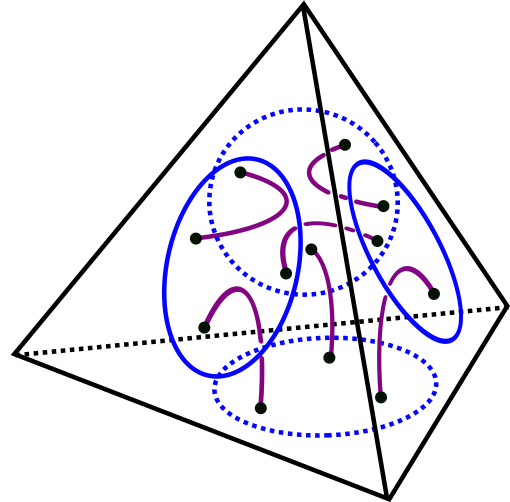}};
        \end{scope}
\end{tikzpicture}
\caption{{\bf Left:} Dihedral angles for hyperbolic polyhedron homeomorphic to $\Delta \sss (J_0 \cap \Delta)$. {\bf Right:} The tangle $J_0 \subset T_0$.}
\label{fig:tangle}
\end{figure}

First, we fix a ``model" tetrahedron $T_0$, and consider the tangle $J_0 \subset T_0$ shown on right of \Cref{fig:tangle}.  Here $J_0 \subset T_0$ consists of six properly embedded arcs, one connecting every pair of faces, and four circles, one contained in each face.  In particular, we can index the components of $J_0$ by the set of faces (for circle components) and edges (for arc components, index an arc by the edge shared by the faces it meets).  In \cite{lw} we showed that $T_0 \sss J_0$, viewed as a manifold with corners along the edges, admits a hyperbolic structure for which the boundary is totally geodesic so that the dihedral angles along each edge are equal to $\frac{\pi}2$. The tangle and hyperbolic structure are constructed to be invariant by the full symmetry group $\Tet \cong S_4$ of $T_0$.  Specifically, we take a fundamental domain $\Delta$ for the action of $\Tet$ on $T_0$, which is one tetrahedron of the first barycentric subdivision of $T_0$.  By inspection, $\Delta$ meets the tangle $J_0$ in exactly two arcs as indicated on the left in \Cref{fig:tangle}.  We then explicitly constructed a partially ideal hyperbolic polyhedron which is homeomorphic to $\Delta \sss J_0$ so that the dihedral angles are as indicated in the figure.  The hyperbolic structure on $T_0 \sss J_0$ is such that the restriction to $\Delta \sss J_0$ (or any image by $\Tet$) is isometric to this partially ideal hyperbolic polyhedron.  The link $L \subset M$ associated to the triangulation $\tau$ is defined by the property that for every tetrahedron $T$, there is a homeomorphism of pairs $(T,T \cap L) \cong (T_0,J_0)$.  The locally $\CAT(-1)$ metric on $M_L$ has $T \sss L$ isometric to the fixed hyperbolic metric on $T_0 \sss J_0$, for every tetrahedron $T$ of $\tau$. 

\begin{figure}[h] 
    \centering
    \begin{tikzpicture}
    \node at (0,0) {\includegraphics[width=\linewidth]{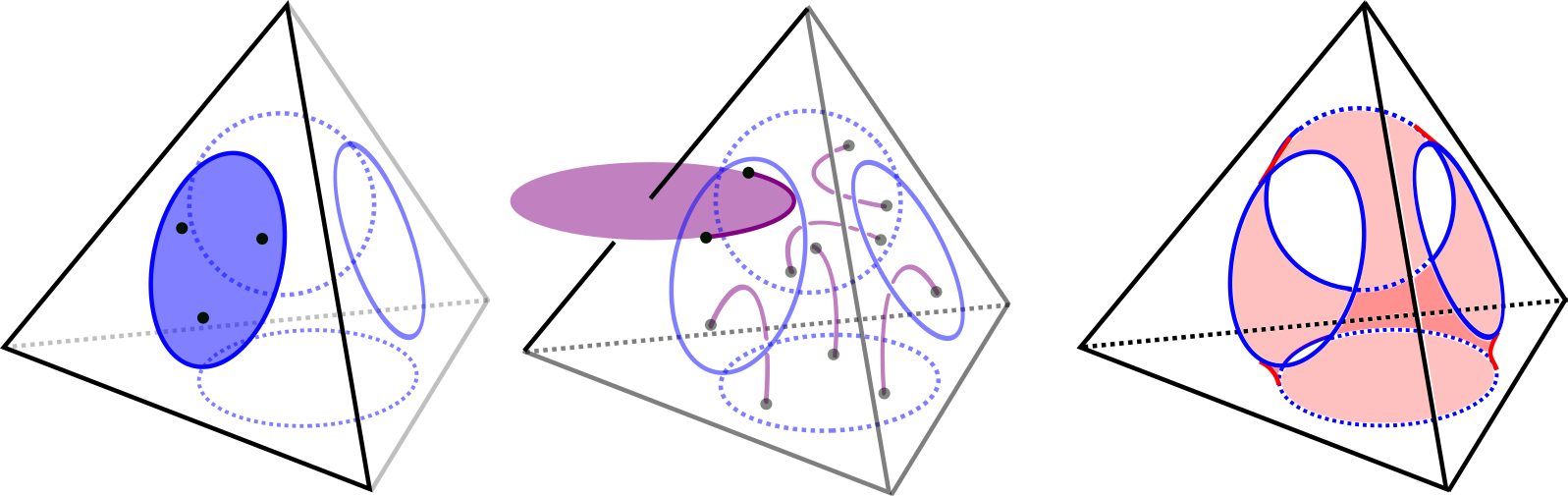}};
    \node at (-4.8,-.25) {$\Sigma_f$};
    \node at (-5.25,-0.9) {\blue{$f$}};
    \node at (-1.4,.4) {$\Sigma_e$};
    \node at (-.9,1.) {$e$};
    \node at (4.6,.9) {$\Sigma_T$};
    \node at (5.8,1) {$T$};
    \end{tikzpicture}
    \caption{\textbf{Left:} A face $f$ and its associated surface $\Sigma_f$. \textbf{Middle:} An edge $e$, and its associated surface $\Sigma_e$. \textbf{Right:} The four-holed sphere $\Sigma_T$ associated with a tetrahedron $T$. }
    \label{fig:surfaces}
\end{figure}

Suppose $L \subset M$ is a link associated to a Cooper-Thurston triangulation, $\tau$, as above.  There are three families of surfaces in $M$, pictured in \Cref{fig:surfaces}, described as follows.  Let $\E_\tau$, $\F_\tau$, and $\T_\tau$ denote the set of edges, faces, and tetrahedra of $\tau$.
    \begin{enumerate}
    \item For every $f \in \F_\tau$ of $\tau$, the \emph{face surface/disk} associated with $f$, denoted $\Sigma_f$, is the disk bounded by the component of $L$ embedded in $f$.
    \item For every $e \in \E_\tau$, there is a component of $L$ that encircles $e$, built from a subset of the arcs of $L$ intersected with the tetrahedra containing $e$. The \emph{edge surface/disk} associated with $e$, denoted $\Sigma_e$, is the disk bounded by this link component.
    \item For every $T \in \T_\tau$, let $f_1,\ldots,f_4$ be the faces of $T$.  The \emph{tetrahedon surface} associated with $T$, denoted $\Sigma_T$, is obtained by ``pushing" the interior of the closure of the complement of the face surfaces in the boundary, $\overline{\partial T \sss \bigcup\limits_{i = 1}^4 \Sigma_{f_i}}$, into the interior of $T$.
    \end{enumerate}
For each edge, face, or tetrahedron $x$, we write $\Sigma^\circ_x$ for the intersection of $\Sigma_x$ with $M_L$. Each $\Sigma_x^\circ$ is homeomorphic to a sphere minus some number of points (or equivalently closed disks) and are pictured in~\Cref{fig:surfaces}. 

\begin{remark} \label{rem:disk v sphere} In either $M_L$ (the locally $\CAT(-1)$ manifold) or $M_L^h$ (the complete hyperbolic manifold) all punctures and (removed) boundary components of face, edge, and tetrahedron surfaces are cusps. Thus, although it is convenient to view the edge and face surfaces as punctured disks from a topological perspective (and the tetrahedron surface as having four ``missing" boundary components), from a geometric point of view, they are more naturally all punctured spheres. In this paper, we will move freely between these two perspectives depending on which is more descriptive. 
\end{remark}

Because every tetrahedron $T$ of $\tau$ has $T \sss L$ isometric to $T_0 \sss J_0$ in the locally $\CAT(-1)$ metric, as explained in \cite{lw}, these surfaces are all totally geodesic.  We recall the reasoning and point out some additional features that we will need.

First, the face surfaces and edge surfaces are constructed by gluing together polygons in the faces of the ideal polyhedron $\Delta \sss J_0$ along sides that have dihedral angle $\frac{\pi}2$ in $\Delta$, and so one can see explicitly that they are totally geodesic.  The number of faces of $\Delta \sss J_0$ needed to construct an edge surface in $M_L$ is precisely twice the number of tetrahedra that the edge meets, and because $L$ comes from a Cooper-Thurston triangulation, these surfaces are built from 8, 12, 16, or 20 polygons.  Additionally, in $T_0 \sss J_0$, the isometric $\Tet$ action implies every face surface is isometric to a single model face surface we denote $\Sigma_F^\circ$.  This is a hyperbolic surface, and the stabilizer of a face surface via the isometric action of $\Tet$ shows that it admits a three-fold rotational symmetry, with quotient a $(3,\infty,\infty)$ triangle orbifold (topologically a sphere with two punctures and one orbifold point of order $3$).  In particular, every face surface in $M_L$ admits such a symmetry.  The group $\Tet$ also leaves the tetrahedral surface invariant, and quotienting this surface by the orientation preserving subgroup yields a triangle orbifold. Theorem 2.1 of \cite{as} implies that the tetrahedron surface is also totally geodesic (to apply the theorem, we can treat the surfaces in the boundary as reflectors).  See \cite{lw} for details.

\subsection{An auxiliary tangle}
\label{subsection:aux tangle}

Next we describe an auxiliary tangle in the closed $3$--ball, $J_{\aux} \subset \bar{\mathbb B}^3$, as well as its geometry.  Given a Cooper-Thurston triangulation $\tau$ and associated link $L \subset M$, we will use copies of this tangle to ``join together" the components of $L$ that intersect the face surfaces to produce a new link $L' \subset M$. 

\begin{figure}[h]
\begin{center}
\begin{tikzpicture}
\node at (0,0) {\includegraphics[width=9cm]{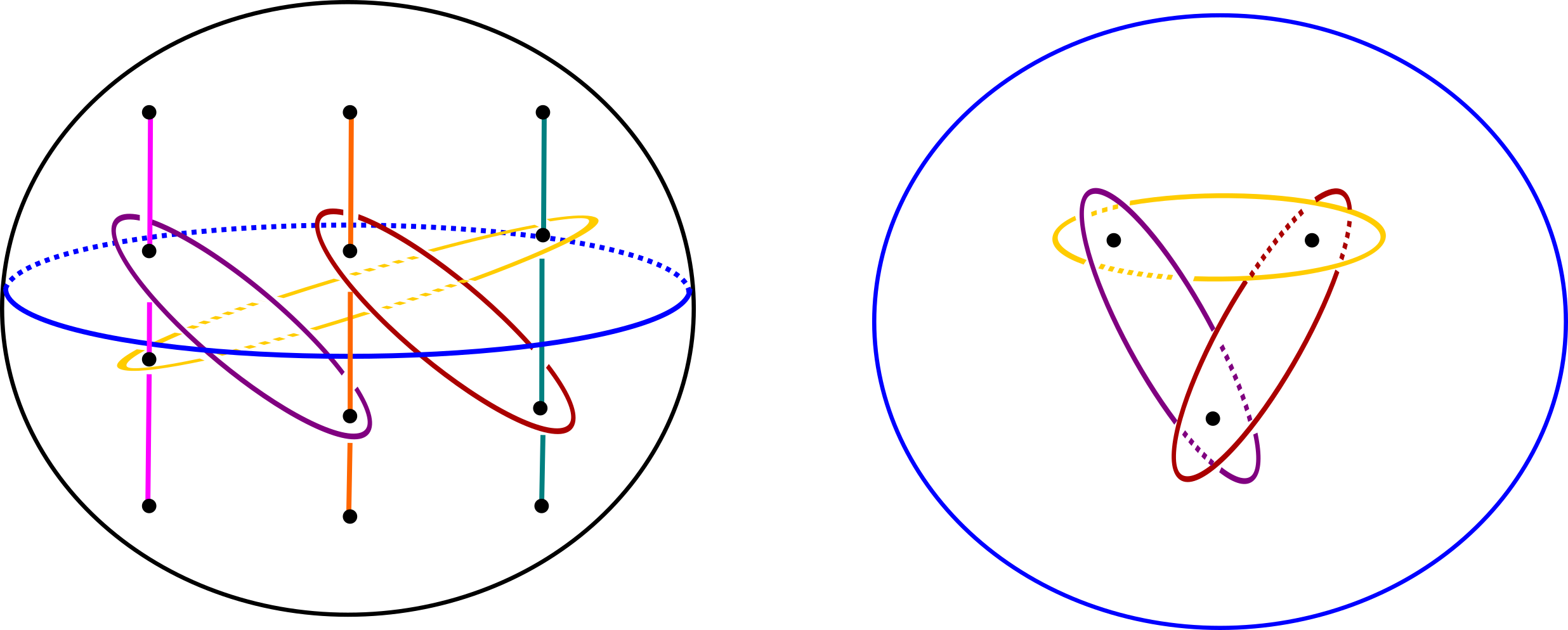}};
\end{tikzpicture}
\caption{Two views of the tangle $J_{\aux} \subset \bar{\mathbb B}^3$, with the view ``from above" on the right, where the arc components are visible only as ``dots".} \label{fig:Aux-tangle-1-new}
\end{center}
\end{figure}

The tangle $J_{\aux}$ is shown in \Cref{fig:Aux-tangle-1-new}.  It consists of three properly embedded arcs, and four simple closed curves.  One of the simple closed curves is on the boundary of $\bar{\mathbb B}^3$, and is illustrated as the equator in the figure.  The other three, which we refer to as the {\em internal components} of $J_{\aux}$, encircle the arcs as shown.   The picture of $J_{\aux}$ shown on the right more clearly illustrates a $3$--fold rotational symmetry.  The complement $\bar{\mathbb B}^3\sss J_{\aux}$ has two boundary components which are each thrice punctured disks (or four-times punctured spheres; c.f.~Remark~\ref{rem:disk v sphere}).

We first make precise some terminology that will be used throughout the rest of this section.  Suppose $\Sigma \subset N$ is a totally geodesic surface in a complete, singular hyperbolic, locally $\CAT(-1)$ manifold $N$.
\begin{defn} \label{def:cut open}
    To {\em cut a manifold $N$ open along a surface $\Sigma$} means to remove $\Sigma$, and then take the path-metric completion of the result.  We write $N/\!/\Sigma$ to denote {\em $N$ cut along $\Sigma$}.
\end{defn} 
The manifold $N/\!/\Sigma$ has two new, totally geodesic boundary components, which may be isometrically glued back together in a way that is canonically isometric to the original manifold $N$.  A similar definition makes sense for cutting open any manifold along a codimension $1$ submanifold.

\begin{theorem}\label{thm:connectorhyperbolic}
    The manifold $\bar{\mathbb B}^3 \sss J_{\aux}$ has a hyperbolic structure for which both boundary components and the three twice-punctured disks bounded by the interior components are totally geodesic.  Furthermore, the $3$--fold symmetry group acts by isometry, and each boundary component is isometric to $\Sigma_F^\circ$, the model face surface in $T_0 \sss J_0$, in a way that respects the three-fold symmetry in each.
\end{theorem}
The thrice punctured disks have empty boundary as surfaces in $\bar {\mathbb B}^3 \sss J_{\aux}$, so from a geometric point of view, they are more naturally four--times punctured spheres, but this nomenclature of ``punctured disks" is helpful for visualizing the surfaces (c.f. Remark~\ref{rem:disk v sphere}).

\begin{proof}
We consider the auxiliary link $L_{\aux} \subset S^3$, shown in \Cref{fig:auxtangle}. The link $L_{\aux}$ can be seen to be hyperbolic by using the \lstinline{verify_hyperbolicity} function in SnapPy \cite{Hoffman}.  Let $D_0$ be the disk in $S^3 \sss L_{\aux}$ bounded by the blue component, as indicated, and write $D_0^\circ = D_0 \cap S^3 \sss L_{\aux}$.

\begin{figure}[h]
    \centering
    \begin{tikzpicture}
    \node at (0,0) {\includegraphics[width=0.6\linewidth]{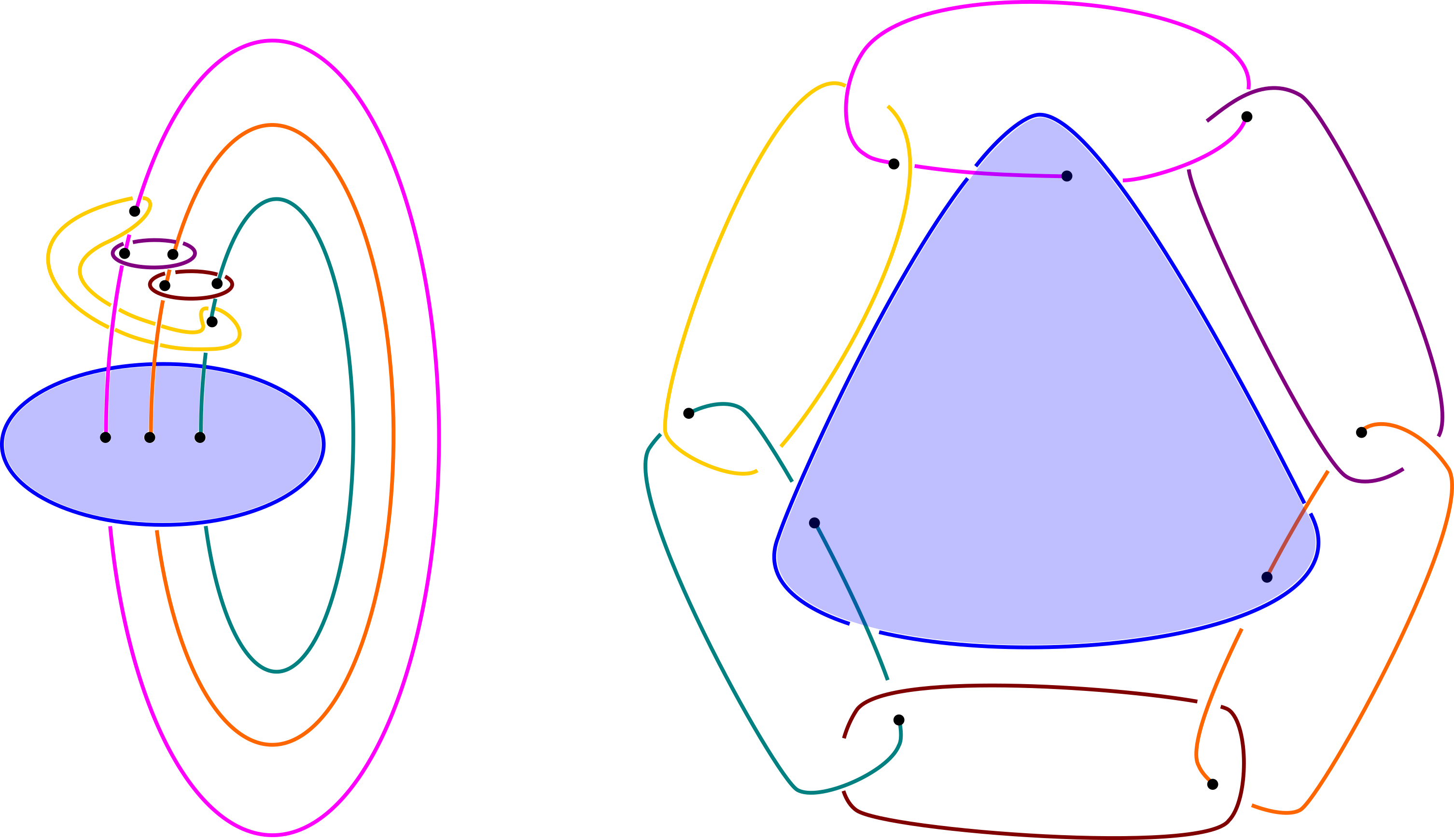}};
    \node at (1,-.5) {$D_0$};
    \node at (-3.25,-0.3) {$D_0$};
\end{tikzpicture}
    \caption{Two views of the link $L_{\aux}$. Note that $(L_{\aux} \sss \partial D_0) \subset S^3$ is a $6$ component chain link, $C_6 \subset S^3$.}
    \label{fig:auxtangle}
\end{figure}

The peripheral loops of $D_0^\circ$ are homotopic to distinct cusps of $S^3 \sss L_{\aux}$, and it follows that all non-null homotopic simple closed curves in $D_0^\circ$ are non-null homotopic in $S^3 \sss L_{\aux}$; hence, $D_0^\circ$ is essential.
Cutting open along $D_0^\circ$ yields $\bar{\mathbb B}^3 \sss J_{\aux}$, and hence hyperbolicity of $\bar{\mathbb B}^3 \sss J_{\aux}$ follows.

To see that the boundary components of $\bar{\mathbb B}^3 \sss J_{\aux}$ are totally geodesic with respect to the given hyperbolic structure, it suffices to show that $D_0^\circ$ is isotopic to a totally geodesic surface. For this, we note that the $3$--fold symmetry of $\bar{\mathbb B}^3 \sss J_{\aux}$ extends to a $3$--fold symmetry of $S^3 \sss L_{\aux}$, which can be realized as an isometry by Mostow's Rigidity Theorem.  Moreover, $D_0^\circ$ is invariant, and the quotient is a $(3,\infty,\infty)$--triangle $2$--orbifold $Q$ in the quotient $\mathcal O$ of $S^3 \sss L_{\aux}$.  Since $D_0^\circ$ is essential in $S^3 \sss L_{\aux}$, $Q$ is essential in $\mathcal O$.
By Theorem 2.1 in \cite{as}, $\mathcal O$ is isotopic to a totally geodesic $2$--orbifold $Q'$.  This isotopy can be lifted to an isotopy of $D_0^\circ$ to a totally geodesic surface ${D_0^\circ}'$.  Furthermore, because the hyperbolic structure on a triangle orbifold is unique by Corollary 13.3.7 of \cite{Thurstonnotes}, it follows that ${D_0^\circ}'$ is isometric to $\Sigma^\circ_F$, as one can lift the isometry between $Q'$ to the quotient of $\Sigma^\circ_F$ by its $3$--fold symmetry group.  

Each of the three internal components of $J_{\aux}$ bound disks, $D_{12}$, $D_{23}$, and $D_{31}$ in $\bar{\mathbb B}^3$ with boundary compoents in $J_{\aux}$.  The intersection of $D_{ij}$ with $\bar {\mathbb B}^3 \sss J_{\aux}$ is a thrice punctured sphere, $D_{ij}^\circ$ in the hyperbolic structure in $\bar{\mathbb B}^3 \sss J_{\aux}$, which are also thrice-punctured spheres in $S^3 \sss L_{\aux}$ (c.f.~Remark~\ref{rem:disk v sphere}).  By \cite[Theorem~2.1]{adams}, after an isotopy we can assume these surfaces are totally geodesic.  Since totally geodesic surfaces have the minimal number of components of transverse intersections in their isotopy classes, and since the $D_{ij}^\circ$ can be realized to be pairwise disjoint after an isotopy, and disjoint from $D_0^\circ$, the same is true of their totally geodesic representatives.  In particular, each $D_{ij}^\circ$ can be assumed to be totally geodesic in $\bar{\mathbb B}^3 \sss J_{\aux}$.
\end{proof}

\subsection{New Links} \label{S:new links}
We now return to the construction of a new link associated with a Cooper-Thurston triangulation $\tau$ of $M$.  We start with the associated link $L = L(\tau) \subset M$ described in \S\ref{subsection:fillinglinkrecap}.  For each face $f$ of $\tau$, we consider two copies of the disk $\Sigma_f$ in $M$ with the same boundary, and slightly push their interiors into the two tetrahedra sharing the face $f$.  The union of each such pair bounds a ball $\bar{\mathbb B}^f$, and $(\bar{\mathbb B}^f,\bar{\mathbb B}^f \cap L)$ is homeomorphic to the tangle $(\bar{\mathbb B}^3,J_{\aux}')$, where $J_{\aux}'$ is the union of components of $J_{\aux}$ {\em not including the three internal components}.  The new link $L'' = L''(\tau)$ is obtained from $L$ by adding three components for each face $f$, so that $(\bar{\mathbb B}^f,\bar{\mathbb B}^f \cap L'') \cong (\bar{\mathbb B}^3,J_{\aux})$.  We note that $L''$ can be viewed as a kind of {\em augmentation} of $L$ as in Adams \cite{Adams-augment}, Lackenby \cite{Lackenby-augment}, Purcell \cite{Purcell-augment}, and many others, and we will utilize this perspective.

Another perspective is that $M_{L''} = M \sss L''$ is homeomorphic to the $3$--manifold obtained by first cutting $M_L$ along the union of the face surfaces (see \Cref{def:cut open}), producing a $3$--manifold with two totally geodesic boundary components  for each face $f \in \mathcal{F}_\tau$, each isometric to $\Sigma_F^\circ$, then gluing in a copy of $\bar{\mathbb B}^3 \sss J_{\aux}$ by isometry along the faces. 
In particular, the link complement $M_{L''}$ contains {\em two} face surfaces for every face $f$, which we denote $\Sigma_{f\pm}^\circ$, in addition to the tetrahedra surfaces $\Sigma_T^\circ$ for each tetrahedron $T \in \mathcal{T}_\tau$.  For every face $f$ of $\tau$, there are also three twice-punctured disks $D_{f12}^\circ$, $D_{f23}^\circ$, and $D_{f31}^\circ$, coming from the punctured disks $D_{12}^\circ$, $D_{23}^\circ$, $D_{31}^\circ$ in $\bar{\mathbb B}^3 \sss J_{\aux}$.  The following is immediate from the discussion above.

\begin{theorem}
    Suppose $L'' = L''(\tau)$ is the link associated to a Cooper-Thurston triangulation, $\tau$, as above.  Then $M_{L''}$ admits a locally $\CAT(-1)$ metric so that for all tetrahedra $T \in \T_\tau$, faces $f \in \F_\tau$, and $ij = 12,23,31$, the surfaces $\Sigma_T^\circ$, $\Sigma_{f\pm}^\circ$, and $D_{fij}^\circ$ are totally geodesic.   Furthermore, these totally geodesic surfaces are pairwise disjoint in $M_{L''}$. \qed
\end{theorem}

\begin{figure}[h]
\begin{center}
\begin{tikzpicture}
\begin{scope}[scale = .8]
    \node at (0,0) {\includegraphics[width=12cm]{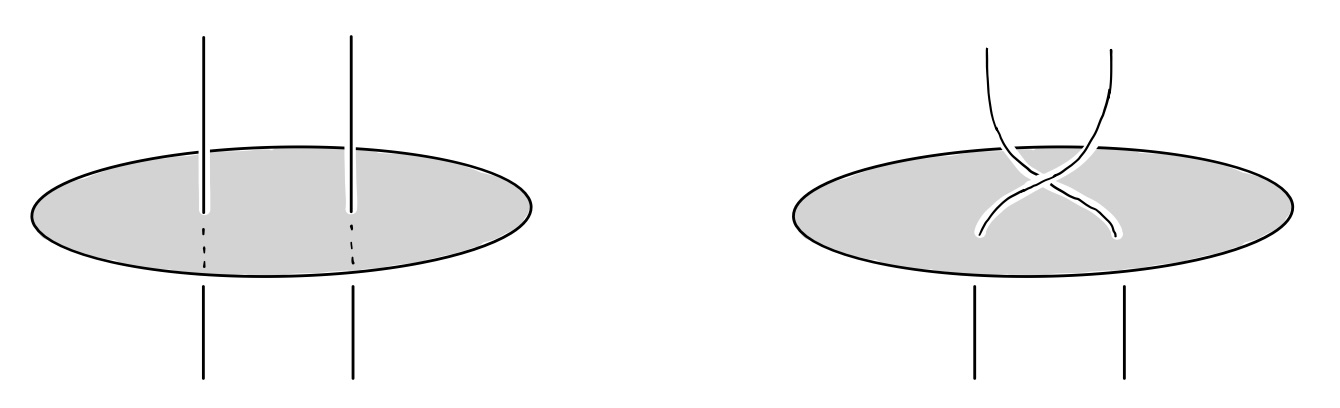}};
    \node at (-5.9,-.1) {$D_{fij}$};
    \draw[->,>=stealth] (-.8,0) -- (.7,0);
    \node at (0,-.5) {\small{half twist}};
    \node at (0,-1) {\small{in $D_{fij}$}};
\end{scope}
\end{tikzpicture}
\caption{Twisting in a disk.}\label{fig:disk twist}
\end{center}
\end{figure}

We now describe the construction of a third link $L' = L'(\tau) \subset M$ associated with $\tau$, from which the knots used to prove the main theorem will be obtained by Dehn surgery.  For this construction, for any face $f \in \F_\tau$ and any $ij$, we let $D_{fij}\subset M$ denote the closure of the punctured disk $D_{fij}^\circ$; we assume $D_{fij}$ is a disk whose boundary is a component of $L''$, and whose interior meets $L''$ in two points.  We can modify $L''$ locally in a neighborhood of $D_{fij}$, by a {\em half twist} along $D_{fij}$, which has the effect of replacing the local picture of $L''$ near $D_{fij}$ on the left of \Cref{fig:disk twist} with the one on the right of that same figure, and keeping the link unchanged elsewhere.

Since the set of disks $\{D_{fij} \}_{ij,f}$ are pairwise disjoint, we can independently perform half-twists on any collection of these disks to produce a new link which depends (up to isotopy) only on the choices of disks in which we perform a half-twist, and not the order in which we performed the twists.  Let $\mathcal D$ denote some set of such disks and $L' \subset M$ the link obtained by half-twisting in those disks.  Observe that we can obtain $M_{L'}$ from $M_{L''}$ by cutting $M_{L''}$ open along the set of punctured disks $D_{fij}^\circ$ whose closures are in $\mathcal D$, and then regluing by precomposing ``the identity" gluing with the half twist in the disk. We consider all $D_{fij}^\circ$, $\Sigma_{f\pm}^\circ$, and $\Sigma_T^\circ$ as surfaces in $M_{L'}$ (and their closures as surfaces in $M$).  The following lemma specifies a choice of $\mathcal{D}$ from which we will construct a specific link $L'$. 

\begin{lemma} \label{lem:final new link}
    There is a choice of disks $\mathcal D$ as above, so that the resulting link $L' \subset M$ consists of the following union,
    \[ L' = K' \sqcup \bigsqcup_{f \in \F_\tau} \left( \partial \Sigma_{f+} \sqcup \big( \bigsqcup_{ij} \partial D_{fij}\big)\right),\]
    where $K'$ is a single component.

    The link complement $M_{L'}$ admits a complete, singular hyperbolic, locally $\CAT(-1)$ metric isometric to $M_{L''}$ outside the union of all the punctured disks ${D_f^\circ}_{ij}$ (in both manifolds).
    In particular, the totally geodesic surfaces $\Sigma_T^\circ$ for $T \in \T_\tau$ and $\Sigma_{f\pm}^\circ$ for $f \in \F_\tau$ in $M_{L''}$ define totally geodesic surfaces of the same name in $M_{L'}$.
\end{lemma}
The last disjoint union is over the indices $ij \in \{12,23,31\}$.  We note that $\partial \Sigma_{f+} = \partial \Sigma_{f-}$.  Assuming the conclusion of the lemma, we observe that for every disk $D$ of the form $\Sigma_{f\pm}$ or $D_{fij}$, the intersection of $L'$ with the interior of $D$ is precisely $K' \cap D$.  This observation is the basis of the proof of the lemma. 

\begin{proof}
Given a choice of disks $\mathcal D$ and associated link $L'$, we say that a disk $D_{fij}$ is {\em good} for $L'$ if the two intersection points with the interior of $D_{fij}$ are from the same component of $L'$.  For example, observe that $L''$ is the link associated to $\mathcal D = \emptyset$, and that by inspection, no disk $D_{fij}$ is good for $L''$.

We first prove the following.
\begin{claim}
    There is a collection $\mathcal D$ such that {\em every} disk $D_{fij}$ is good for the associated link $L'$.
\end{claim}
\begin{proof}[Proof of claim]
First, enumerate all the disks, over all faces $f \in \tau$ and all $ij \in \{12,23,31\}$, writing the list as $D_1,\ldots,D_m$.  We will define the collection $\mathcal D$ as a non-decreasing union
\[ \mathcal D_1 \subseteq \mathcal D_2 \subseteq \mathcal D_3 \subseteq \cdots \mathcal D = \bigcup_{k=1}^m \mathcal D_k.\]
The sets $\mathcal D_k$ will be defined recursively, having the property that
\[ \mathcal D_k \subset \{ D_1,\ldots,D_k\},\]
{\em and} so that $\mathcal D_k$ is obtained from $\mathcal D_{k+1}$ by either adding $\overline D_k$ to $\mathcal D_{k-1}$ or not.  For each $k$, we let $L'_k$ be the link associated to $\mathcal D_k$.  The properties of the sequence $\mathcal D_k$ described above implies that $L'_k$ is obtained from $L'_{k-1}$ by either twisting in $D_k$ or not, depending on whether it is included in $\mathcal D_k$ or not.

Now we define $\mathcal D_1 = \{D_1\}$, and since $D_1$ is not good for $L''$, it {\em is} good for $\mathcal D_1$: the half-twist has the effect of combining the two components of $L''$ meeting the interior of $D_1$ into a single component of $L_1'$ (see \Cref{fig:disk twist}).

Next, having defined $\mathcal D_{k-1}$ for some $k >1$, we define 
\[ \mathcal D_k = \left\{ \begin{array}{ll} \mathcal D_{k-1} & \mbox{if $D_k$ is good for $L_{k-1}'$}\\
\mathcal D_{k-1} \cup \{D_k\} & \mbox{otherwise} \end{array} \right.\]
We observe that, by induction, for all $1 \leq j \leq k$, $D_j$ is good for $L_k$, since a half twist in a disk that isn't good makes it good (and doesn't destroy goodness for any other disks). The set $\mathcal{D}$ thus satisfies the property in the claim. 

The second half of the statement of Lemma \ref{lem:final new link} follows immediately from Corollary 4.4 of ~\cite{adams}. Indeed, we obtain $M_{L'}$ from $M_{L''}$ by cutting along, applying half-twists to, and regluing totally geodesic thrice-punctured disks. Since each half-twist is an isometry for the hyperbolic metric on the punctured disk, the resulting complete, singular hyperbolic, locally $\CAT(-1)$ manifold is isometric to $M_{L''}$ in the complement of those punctured disks along which we performed the half twists.
\end{proof}

For each face $f$, there are three arcs of $L'$ that connect the interior of $\Sigma_{f-}$ to the interior of $\Sigma_{f+}$, and each one of these arcs intersects two of the disks $D_{fij}$.  Since each disk is good, it follows that all three of these arcs are contained in the same component of $L'$.  From this we see that for any tetrahedron $T \in \T_\tau$, if $f_1,\ldots,f_4$ are the faces, then all arcs of $L'$ that meet the interiors of the disks $\Sigma_{f_i\pm}$ ($i=1,2,3,4$) belong to the same component of $L'$.  Connectivity of $M$ implies that there is a single component $K'$ of $L'$ with the property that for any disk $D_{fij}$ or $ \Sigma_{f\pm}$, the intersection of $L'$ with its interior is precisely the intersection of $K'$ with its interior.  Thus, $L'$ is the union of $K'$ and the boundaries of these disks.
\end{proof}

\begin{defn} \label{def:main link}
    Given a Cooper-Thurston triangulation $\tau$ of a closed orientable $3$--manifold, $M$, we write $L' = L'(\tau)$ for the link coming from some choice of $\mathcal D$ as in Lemma~\ref{lem:final new link}, and refer to this as a {\em fully twisted} link associated to $\tau$.  We will refer to the component $K' \subset L'$ specified in the lemma as the {\em focal component} and all other components as {\em disk components}.
\end{defn}

Given a fully twisted link $L'$ associated to $\tau$, it will be useful to have some language to describe how $L'$ sits inside $M$, and its relation to $L$ and the triangulation $\tau$.  Recall that for each face $f$, $\bar{\mathbb B}^f$ is the ball in $M$ bounded by the union of faces $\Sigma_{f+} \cup  \Sigma_{f-}$.  The collection of balls $\{\bar{\mathbb B}^f\}_f$ over all faces is a set of pairwise disjoint balls.  

For each face $f$, the pair $(\bar{\mathbb B}^f,\bar{\mathbb B}^f \cap L')$ is homeomorphic to $(\bar{\mathbb B}^3,J_{\aux}^f)$, where $J_{\aux}^f$ is obtained from $J_{\aux}$ after possibly applying some half-twists to the disks $D_{ij}$ (namely those specified by the choice of $\mathcal D$).  

\begin{defn} \label{def:connector}
    We write
    \[ \mathbb B^f = \bar{\mathbb B}^f \sss (\bar{\mathbb B}^f \cap L'),\]
    and call this manifold the {\em connector for $f$.}
\end{defn}

It is useful to think about each tetrahedron as having been cut along the interiors of $\Sigma_f$, over its four faces $f$, and pushed off of the interiors of the associated balls $\bar{\mathbb B}^f$.
In particular, $M$ is reconstructed by gluing together the disjoint union of the tetrahedra from $\tau$ and the balls $\{\bar{\mathbb B}^f\}$.  Thus, we think of removing the union of the interiors of all  $\bar{\mathbb B}^f$ as cutting $M$ open along the interiors of the face disks.  

Consequently, the way $L'$ meets the complement of the interiors of $\bar{\mathbb B}^f$ in $M$ is precisely the way $L$ meets $M$ cut along the interiors of the face disks.  More precisely, let us write $\check M$ to denote the complement of the unions of the interiors of $\bar{\mathbb B}^f$ in $M$,
\[ \check M = M \sss \bigcup_f \mbox{int}(\bar{\mathbb B}^f),\]
and $N$ to denote $M$ cut along the interiors of the face surfaces, \[N = M /\!/ \bigcup_f \mbox{int}\big( \Sigma_f\big)\] (as in \Cref{def:cut open}), and then
\begin{equation} \label{eq:L and L'}
    (\check M,\check M \cap L') \cong (N, N\cap L).
\end{equation}

In particular, we can view $M_{L'}$ as obtained from $M_{L}$ by cutting along each $\Sigma_f^\circ$, and gluing in the $f$--connector, $\mathbb B^f$,
\begin{equation} \label{eq:L and L' complement} M_{L'} \cong \big( M_L/\!/ \bigcup_f \Sigma_f^\circ \big) \cup \bigcup_f \mathbb B^f \cong N_L \cup \bigcup_f \mathbb B^f,
\end{equation}
where $N_L = N \sss (N \cap L)$.




\subsection{Constructing knots and the theorems revisited}
\label{sec:constructingknots}

Given a Cooper--Thurston triangulation $\tau$ in a closed $3$--manifold $M$, and an integer $n$, we define a knot $K(n,\tau) \subset M$ as follows.  Let $L'$ be the fully twisted link associated to $\tau$ and $K'$ its focal component.
The disk components of $L'$ (that is, the components in $L' \sss K'$) bound disks with pairwise disjoint interiors 
\[
\bigcup_{f,ij} \Sigma_{f+} \cup \Sigma_{f-} \cup D_{fij}.
\]
We assume that the peripheral tori of $M_{L'}$ are framed with basis for which the longitude is defined by these disks (and the meridian is the one coming from the solid tori $\mathcal N(L')$).  

Now consider the link $L' \sss K' \subset M$ consisting of {\em only the disk components of $L'$}.  For any integer $n$, we apply $(1,n)$--Dehn surgery on {\em every} component of this link.  Since the disks bounded by the disk components have pairwise  disjoint interiors, the resulting manifold is homeomorphic to $M$.   The image of  the focal component $K'$ by this homeomorphism is a knot which we denote $K(n) = K(n,\tau)$.  In fact, we can view $K(n)$ as obtained from $K'$ by adding $n$ full twists at each disk (i.e.~$2n$ half twists as in \Cref{fig:disk twist}); see \cite{Adams-augment,Lackenby-augment,Purcell-augment} for this kind of construction.
The homeomorphism sends the disk components to the core curves of the filling solid tori, which results in a link  that contains $K(n)$ as one of components and whose components all bound disks. 

\begin{remark}
    We can view $M_{K(n)} = M \sss K(n)$ as obtained by applying a $(1,n)$--Dehn filling on each cusp of $M_{L'}$ coming from a disk component of $L'$, and leaving the cusp corresponding to the focal component unfilled; that is
    \[ M_{K(n)} \cong M_{L'}((1,n),\ldots,(1,n),\infty).\]
\end{remark}

We can now state a more precise version of the theorems from the introduction.  The first of these immediately implies \Cref{thm:longessentialsys}.

\begin{theorem}\label{thm:main1specific}
    Suppose $M$ is any closed, orientabled $3$--manifold with $\pi_1M \neq \{\ast\}$ and $R > 0$.  Then there exists a Cooper--Thurston triangulation $\tau$ so that for all $n$ sufficiently large, the associated knot $K(n,\tau)$ has
    \[ \esssys(K(n,\tau)) > R.\]
\end{theorem}

A map $f \colon \Gamma \to M$ has {\em full rank} if $\rank(f_*(\pi_1(\Gamma))) = \rank(\pi_1(\Gamma))$. A knot $K \subset M$ is {\em full rank--$m$ filling} if for every $f \colon \Gamma \to M_{K}$ such that $\rank(\pi_1(\Gamma)) \leq m$ and $i \circ f \colon \Gamma \to M$ is full rank, $f_*$ is injective. The next theorem implies \Cref{T:filling knot} by setting $m = \rank(\pi_1(M))$.

\begin{theorem} \label{thm:main2specific}
    Suppose $M$ is any closed, orientable $3$--manifold and $m > 0$.  Then there exists a Cooper--Thurston triangulation $\tau$ so that for all $n$ sufficiently large, the associated knot $K(n,\tau)$ is full rank--$m$ filling.
\end{theorem}

\section{Hyperbolic knots with large essential geodesic systole} \label{sec:longessentialsystole}

From now on, the paper diverges in two different directions. In this section, we prove \Cref{thm:main1specific} using \emph{hyperbolic} structures on $K(n, \tau)$. In the section that follows, we prove \Cref{thm:main2specific} using a locally $\CAT(-\tfrac12)$ structure rather than a hyperbolic structure. 

The crux of this section is the proof that the link $L'=L'(\tau)$ constructed from $\tau$ in Lemma \ref{lem:final new link} has large essential geodesic systole, provided the combinatorial systole is sufficiently large.

\begin{proposition}\label{prop:predehnfillinglongessentialsys} For $R_0 > 0$, there exists $r > 0$ so that if $\tau$ is a Cooper-Thurston triangulation of $M$ with $\combsys(\tau) > r$, then the hyperbolic link $L' = L'(\tau) \subset M$ has
\[ \esssys(L') > R_0.\] 
\end{proposition}
We postpone the proof of this proposition, and first use it to prove \Cref{thm:main1specific}.

\begin{proof}[Proof of \Cref{thm:main1specific}] Fix $R_0 > 0$. By Proposition \ref{prop:predehnfillinglongessentialsys} and Proposition 3.3 of \cite{lw}, one may take a Cooper-Thurston triangulation of $M$ and its associated link $L'$ with the property that $\esssys(L') > 2R_0$. Given $n \in \mathbb{N}$, recall that \[M_{K(n)} \cong M_{L'}\left((1,n), \ldots, (1,n),\infty\right)\] denotes the manifold obtained by applying the $(1, n)$-Dehn fillings described in \Cref{sec:constructingknots}. By \Cref{thm:dehnfilling}, since $1 + n^2 \rightarrow \infty$ as $n \rightarrow \infty$, $M_{K(n)}$ is hyperbolic, and converges with respect to the geometric topology to the hyperbolic manifold $M_{L'}$. More precisely, there exists a sequence $\epsilon_n \rightarrow 0$ with the property that there is a $(1 + \epsilon_n)$--biLipschitz diffeomorphism between the $\epsilon_n$-thick parts of $M_{K(n)}$ and $M_{L'}$. Additionally, observe that for some $\epsilon' < \epsilon_0$ (the Margulis constant) and $n$ sufficiently large, $(M_{K(n)})_{(0, \epsilon')}$ consists of a cusp neighborhood corresponding to the knot $K(n)$ and the union of solid tori which are the filling solid tori (up to isotopy).

For every $n \in \mathbb{N}$, let $\gamma_n$ denote a geodesic whose length attains the essential geodesic systole of $K(n)$. We claim that there exists some $N \in \mathbb{N}$ such that \[\esssys\left(M_{K(N)}\right) = \ell_{M_{K(N)}}(\gamma_N) > R.\] We split our proof into two cases.

First, suppose that for all $\epsilon > 0$, there exists $n > 0$ such that $\gamma_n$ is not contained in $\left(M_{K(n)}\right)_{[\epsilon, \infty)}$. For each $j$ with $\tfrac1j < \epsilon'$, choose some $n_j$ such that $\gamma_{n_j}$ is not contained in $(M_{K(n_j)})_{[1/j, \infty)}$. We first claim that for each $j$, $\gamma_{n_j}$ cannot be contained in $(M_{K(n_j)})_{(0, \epsilon_0)}$. Indeed, if $\gamma_{n_j} \subset (M_{K(n_j)})_{(0, \epsilon_0)}$, then $\gamma_{n_j}$ is contained in either a cusp neighborhood corresponding to $K(n_j)$ or a Margulis tube, which would necessarily be a filling solid torus (because it would contain an $\tfrac1j$--thin part, and $\tfrac1j < \epsilon'$). If $\gamma_{n_j}$ is contained in a cusp neighborhood, then it cannot be geodesic, a contradiction.  On the other hand, if $\gamma_{n_j}$ is included in a filling solid torus $V$, then since the inclusion of $V$ into $M_{K(n_j)}$ is nullhomotopic, $\gamma_{n_j}$ is null-homotopic in $M$, another contradiction. Thus it must be the case that for all $j$, $\gamma_{n_j}$ intersects $(M_{K(n_j)})_{(0, 1/j)}$, but is not completely contained in $(M_{K(n_j)})_{(0,\epsilon_0)}$. By Lemma \ref{lem:thinpart}, $\ell_{M_{K(n_j)}}(\gamma_{n_j}) \geq R(\tfrac1j)$, and $R(\tfrac1j) \rightarrow \infty$ as $j \rightarrow \infty$. From this one can find some $N \in \mathbb{N}$ with the property that $\ell_{M_{K(N)}}(\gamma_N) > R_0$, as desired.  


Second, suppose that there exists some $0 < \epsilon < \epsilon'$ with the property that for all $n > 0$, $\gamma_n \subset (M_{K(n)})_{[\epsilon, \infty)}.$ Choose some $N \in \mathbb{N}$ sufficiently large so that $(M_{K(N)})_{[\epsilon, \infty)} \subset (M_{K(N)})_{[\epsilon_N, \infty)}$ and $\epsilon_N \leq 1$. Note that since $\gamma_N$ is nontrivial in $M$, $\ell_{M_{L'}}(\gamma_N) \geq \esssys(L')$. 
Since there is an $(1 + \epsilon_N)$-biLipschitz map between $(M_{L'})_{[\epsilon_N, \infty)}$ and $(M_{K(N)})_{[\epsilon_N, \infty)}$, we have

\[\ell_{M_{K(N)}}(\gamma_N) \geq \frac{1}{(1 + \epsilon_N)}\ell_{M_{L'}}(\gamma_N) > \frac{2R_0}{1 + \epsilon_N} \geq R_0.\]
\end{proof} 

\subsection{Preliminaries for the proof of Proposition \ref{prop:predehnfillinglongessentialsys}}\label{sec:longcombsys}


We start by recalling some notation and describing the setup for the rest of the section.  For any Cooper-Thurston triangulation $\tau$, we denote its associated links $L \subset M$ and $L' \subset M$, with complements $M_L$ and $M_{L'}$, respectively.  These implicitly denote the manifolds with their locally $\CAT(-1)$ metrics from \cite[Lemma 3.8]{lw} and  Lemma~\ref{lem:final new link}, respectively.  Recall that both $L$ and $L'$ are hyperbolic by \Cref{T:CAT-1 to hyp}, and we write $M_L^h$ and $M_{L'}^h$ for the hyperbolic metrics on $M_L$ and $M_{L'}$, respectively.  

A first observation that will simplify the discussion is the following consequence of Lemma~\ref{lem:thinpart}.  Recall that $\epsilon_0$ is the Margulis constant, and let $R(\epsilon)$ denote the function of $\epsilon$ from Lemma~\ref{lem:thinpart}.
\begin{lemma} \label{lem:where is esssys}
    Given $R_0 > 0$, suppose $0 < \epsilon < \epsilon_0$ is such that $R(\epsilon) > R_0$.  If $\gamma^h$ is a closed hyperbolic geodesic in $M_{L'}^h$ that meets the $\epsilon$--thin part, $(M_{L'}^h)_{(0,\epsilon)}$, then either $\gamma^h$ is contained in $(M_{L'}^h)_{(0,\epsilon)}$, or else has length greater than $R_0$.
\end{lemma}
\begin{proof}
   Suppose $\gamma^h$ meets $(M_{L'}^h)_{(0, \epsilon)}$, but is not contained in it. Since  $(M_{L'}^h)_{(0, \epsilon_0)} \supset (M_{L'}^h)_{(0, \epsilon)}$, and because there are no closed geodesics contained in cusps, $\gamma^h$ must meet the boundary tori of $(M_{L'}^h)_{(0, \epsilon_0)}$ and $(M_{L'}^h)_{(0, \epsilon)}$. Furthermore, since $\gamma^h$ travels between them, it must have length at least $R(\epsilon) > R_0$.
\end{proof}
Given the setup of the lemma, we note that a geodesic in the $\epsilon$--thin part is necessary a multiple of the core curve of a Margulis tube.  Thus, if a $\gamma^h$ is a geodesic in $M_{L'}^h$, 
then one of the following must hold.
\begin{enumerate}
    \item $\gamma^h$ has length greater than $R_0$,
    \item $\gamma^h$ has length less than $\epsilon$, or
    \item $\gamma^h$ is contained in the $\epsilon$--thick part of $M_{L'}^h$.
\end{enumerate}
This is useful to keep in mind for a geodesic, $\gamma^h$, realizing the essential systole, as we prepare for the proof of Proposition~\ref{prop:predehnfillinglongessentialsys}.  In the first case, we are done, and we will be able to fairly easily rule out the second case when $\combsys(L')$ is sufficiently large.  Thus, the majority of the discussion that follows develops ideas to deal with the third case.  


\subsection{Crossing surfaces}
Recall that for a Cooper--Thurston triangulation $\tau$ of $M$, we write $N$ to denotes the manifold obtained by cutting $M$ open along the interiors of the face surfaces. Then $M$ is obtained from $N$ by gluing in the union of the balls $\bar{\mathbb B}^f$ to the exposed sphere boundary components; in particular, $M$ is the union of tetrahedra in $\T_\tau$ and the balls $\bar{\mathbb B}^f$ for $f \in \F_\tau$.  We consider a cell structure for which (the interiors of) the tetrahedra and the balls are $3$--cells, and call this a {\em $\tau$--cell structure}.  This is imprecise as we have not described the $2$--skeleton, but this will be irrelevant (and one can easily be constructed by appropriately subdividing the $2$--skeleton of $\tau$).

The next lemma follows easily from \cite[Lemma~3.2]{lw}.

\begin{lemma} \label{lem:many tetrahedra}
    Given any closed loop $\gamma$ in $M$ which is non-nullhomotopic, $\gamma$ intersects at least $D=\lfloor \frac{\combsys(\tau)}{16}\rfloor$ tetrahedra of a $\tau$--cell structure.  Furthermore, there exists a collection of such tetrahedra $T_1,\ldots,T_D$ with the property that for all $1 \leq i,j \leq D$, $d_{\hat \tau}(T_i,T_j) \geq 8|j-i|$.
\end{lemma}
\begin{proof}
There is a map $M \to M$ obtained by collapsing the balls onto the face disks, mapping each tetrahedron homeomorphically onto the corresponding tetrahedron in $\tau$.  This sends $\gamma$ to a loop $\gamma'$ which is non-null-homotopic.  By \cite[Lemma~3.1]{lw}, $\gamma'$ meets tetrahedra $T_1,\ldots,T_D$ with the required properties.  Hence $\gamma$ meets the corresponding tetrahedra in the $\tau$--cell structure. 
\end{proof}

\begin{figure}[h]
    \centering
    \begin{tikzpicture}
    \node at (0, 0) {\includegraphics[width=0.8\linewidth]{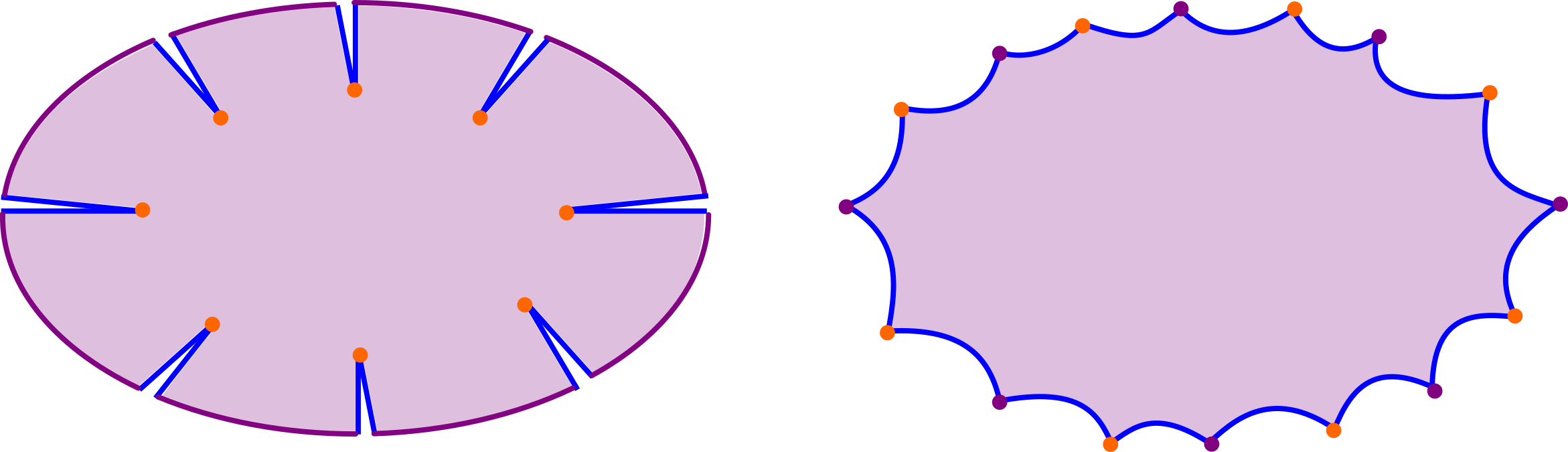}};
    \node at (0, 0) {$\cong$};
    \end{tikzpicture}
    \caption{A cut edge surface $\check{\Sigma}_e^\circ$ may be viewed as an ideal polygon with a uniformly bounded number of sides.}
    \label{fig:sigmacheck}
\end{figure}

Next, recall that by Lemma~\ref{lem:final new link}, there are face and tetrahedron surfaces in $M_{L'}$, denoted
\[ \mathcal S_\tau^\circ(\F,\T) = \{\Sigma_{f\pm}^\circ\}_{f \in \F_\tau} \cup \{\Sigma_T^\circ \}_{T \in \T_\tau}.\]
The edge surfaces from $M_L$ are not defined in $M_{L'}$, but we describe a modification.
Recall that $N_L$ is $M_L$ cut open along the union of its face surfaces, and from \eqref{eq:L and L' complement}, $M_{L'}$ is obtained from $N_L$ by gluing in the connectors $\mathbb B^f$.  For any edge surface $\Sigma_e^\circ$ in $M_{L}$, we can cut $\Sigma_e^\circ$ open and view it inside of $N_L$, and hence inside $M_{L'}$.
More precisely, we write
\[ \check{\Sigma}_e^\circ \subset  N_L \subset M_{L'}\]
to denote the surface obtained by cutting $\Sigma_e^\circ \subset M_L$ along the union of the face surfaces.  This surface is isometric to an ideal hyperbolic $n$--gon, possibly with a cone singularity in the center; see \Cref{fig:sigmacheck}.  Because of the constraints on the number of tetrahedra containing any edge $e$ in a Cooper-Thurston triangulation, $n$ is equal to $8$, $12$, $16$, or $20$ (see \cite[Lemma~3.4]{lw}), and the cone angle at the center is $\tfrac{n}4 \pi$. We call $\check{\Sigma}^\circ_e$ the {\em cut edge surface}, and we write $\mathcal S_\tau^\circ(\E) = \{\check{\Sigma}_e^\circ\}_e$ to denote the set of all cut edge surfaces, over all edges $e$ of $\tau$.

We write $\mathcal S_\tau^\circ = \mathcal S_\tau^\circ(\F,\T) \cup \mathcal S_\tau^\circ(\E)$.  Each of these surfaces is totally geodesic in the locally $\CAT(-1)$ metric in $M_{L'}$.  
Given a tetrahedron $T$, consider the family of surfaces $\mathcal S^\circ_\tau(T) \subset \mathcal S_\tau^\circ$
consisting of the following:
\begin{enumerate}
    \item $\Sigma_{f\pm}^\circ$, where $f$ is a face of $T$,
    \item $\Sigma_{T'}^\circ$, where $T'$ is a tetrahedron and $T' \cap T \neq \emptyset$, and
    \item $\check \Sigma_e^\circ$, where $e$ is an edge of $T$.
\end{enumerate} 
\begin{lemma} \label{lem:tet intersect}
    If $\gamma$ is a closed geodesic in $M_{L'}$ which is non-null-homotopic in $M$, and has non-empty intersection with some tetrahedron $T \in \T_\tau$, then $\gamma$ has non-empty transverse intersection with some surface in $\mathcal S^\circ_\tau(T)$.
\end{lemma}
\begin{remark}
Since we have written $M_{L'}$, rather than $M_{L'}^h$, we note that $\gamma$ is thus assumed to be geodesic in the locally $\CAT(-1)$ metric.
\end{remark}
\begin{proof}
The closures of the surfaces in $\mathcal S_\tau^\circ$ cut $M$ up into simply connected complementary regions in $M$, each of which is one of the following three types: (i) one of the balls $\bar{\mathbb B}^f$, (ii) a region contained in a single tetrahedron, or (iii) a neighborhood of a vertex (thus contained in a collection of tetrahedra that have non-empty pairwise intersection).  Since $\gamma$ is assumed non-null-homotopic in $M$, it cannot be contained in any one of these simply connected regions of $M$.  Likewise, because the (totally geodesic) surfaces in $\mathcal S_\tau^\circ$ are null-homotopic in $M$, $\gamma$ cannot be contained in any of these surfaces, and hence intersects them transversely.  Assuming $\gamma$ intersects $T \in \T_\tau$, it follows that it must intersect one of the simply connected regions of type (ii) contained in $T$, or (iii) which is a neighborhood of a vertex of $T$.  As $\gamma$ must escape this region, it has to intersect one of the surfaces in $\mathcal S_\tau^\circ$, and by definition, such a surface is in $\mathcal S_\tau^\circ(T)$.
\end{proof}

\subsection{Pleating surfaces in $\mathcal S_\tau^\circ$.}
When we pass from the locally $\CAT(-1)$ metric to the hyperbolic metric, the surfaces in $\mathcal S_\tau^\circ$ are no longer necessarily totally geodesic.  Given a geodesic $\gamma$ in $M_{L'}$ that intersects a surface in $\mathcal S_\tau^\circ(\F,\T)$, we will show that the homotopic hyperbolic geodesic, $\gamma^h$ in $M_{L'}^h$, will intersect a pleated representative (Lemma~\ref{lem:meet the tet pleat}).  This is not necessarily true for the edge surfaces, but we still deduce some useful information (Lemma~\ref{lem:near a surface}).

We start by describing the kinds of pleated surfaces we will consider for each $\Sigma^\circ \in \mathcal S_\tau^\circ$.  First, each cut edge surface, $\check \Sigma_e^\circ$, is an ideal hyperbolic polygon in the hyperbolic metric which we denote $\Sigma_e^h$, and we represent it by a pleated surface $\check \Sigma_e^h \to M_{L'}^h$ pleated along an ideal triangulation of $\check \Sigma_e^h$. Given a face surface $\Sigma_{f\pm}^\circ$, there are three sides of cut edge surfaces which meet $\Sigma_{f\pm}^\circ$, and these are pairwise disjoint in $\Sigma_{f\pm}^\circ$.  We extend this union of geodesics to an ideal triangulation, and take the pleated representative we denote $\Sigma_{f\pm}^h \to M_{L'}^h$, pleated along this ideal triangulation.  We realize each tetrahedron surface by any pleated surface, denoted $\Sigma_T^h \to M_{L'}^h$.

\begin{lemma} \label{lem:meet the tet pleat}
    If $\gamma$ is a closed geodesic in $M_{L'}$ that transversely intersects a face or tetrahedron surface $\Sigma^\circ \in \mathcal S_\tau(\F,\T)$, then the hyperbolic geodesic $\gamma^h$ homotopic to $\gamma$ has non-empty intersection with the pleated representative $\Sigma^h \to M_{L'}^h$.
\end{lemma}
\begin{proof}
    By \Cref{T:CAT-1 to hyp}, the identity lifts to a quasi-isometry  $\widetilde M_{L'} \to \widetilde M_{L'}^h \cong \mathbb H^3$ between the universal covers, which extends to a homeomorphism of the Gromov boundaries.  Let $\widetilde \Sigma^\circ$ be a component of the preimage of $\Sigma^\circ$, which is a totally geodesic hyperbolic plane in the $\CAT(-1)$ space $\widetilde M_{L'}$.  The boundary (quasi-)circle at infinity separates the boundary at infinity, $\partial \widetilde M_{L'}$, into two components, and there is an elevation $\widetilde \gamma$ of $\gamma$ with one endpoint in each component.  
    
    Now with respect to the hyperbolic metric, there is a pleated surface $\widetilde \Sigma^h \to \widetilde M_{L'}^h$ covering $\Sigma^h \to M_{L'}^h$ which is a bounded Hausdorff distance from the image of $\widetilde \Sigma^\circ$.  The (quasi-)circle at infinity is the image of the circle which is the boundary at infinity of $\widetilde \Sigma^\circ$.  There is also an elevation $\widetilde \gamma^h$ of $\gamma^h$ which is a bounded Hausdorff distance from the image of $\widetilde \gamma$, and hence has the same endpoints at infinity.  In particular, these are in complementary components of the boundary at infinity of the pleated surface, and hence $\widetilde \gamma^h$ meets that surface.  Consequently, $\gamma^h$ has non-empty intersection with $\Sigma^h \to M_{L'}^h$.
\end{proof}

For each edge $e \in \E_\tau$, each side of the cut edge surface $\check \Sigma_e^\circ$ is contained in a unique face surface $\Sigma_{f\pm}^\circ$.  Observe that the surface $\Sigma_{f+}^\circ$ is such a surface if and only if $\Sigma_{f-}^\circ$ is. We write $\F(e) \subset \F_\tau$ to denote the set of faces $f$ that arise in this way (i.e.~so that $\Sigma_{f\pm}^\circ$ contains a side of $\check \Sigma_e^\circ$), and note that $\F(e)$ is precisely the set of faces that contain $e$.

\begin{lemma} \label{lem:near a surface}
    Suppose $\gamma$ is a closed geodesic in $M_{L'}$ which is non-nullhomotopic in $M$ and has non-empty intersection with a cut edge surface $\check \Sigma_e^\circ$.  Then for some $f \in \F(e)$, the hyperbolic geodesic $\gamma^h$ homotopic to $\gamma$ passes within distance at most $3$ from one of the pleated surfaces $\Sigma_{f+}^h \to M_{L'}^h$ or $\Sigma_{f-}^h \to M_{L'}^h$, for some $f \in \F(e)$.
\end{lemma}
\begin{proof} 
    We consider an elevation of $\gamma$ to a bi-infinite geodesic   $\widetilde \gamma$ in the universal cover $\widetilde M_{L'}$.  Since $\check \Sigma_e^\circ$ is an ideal $n$--gon, which is simply connected, it lifts to the universal cover $\check \Sigma_e^\circ \cong \widetilde \Sigma_e^\circ \subset \widetilde M_{L'}$.  For each $f \in \F(e)$, there is a component of the preimage of $\Sigma_{f+}^\circ$ and $\Sigma_{f-}^\circ$ that meets $\widetilde \Sigma_e^\circ$ along the lifted side.  These are $n$ totally geodesic hyperbolic planes denoted $\{\widetilde \Sigma_{f\pm}^\circ\}_{f \in \F(e)}$ in $\widetilde M_{L'}$.
    The boundaries at infinity of these hyperbolic planes form a set of (quasi-)circles that chain together in a cycle, with the $n$--gon $\widetilde \Sigma_e^\circ$ connecting the consecutive ``tangent" points in the chain.  The $n$ circles bound $n$ disks with pairwise disjoint interiors, $\Omega_1,\ldots,\Omega_n$.  The union of the circles thus divides $\partial \widetilde M_{L'}$ into $n+2$ regions, consisting of the $n$ open disks above, and two remaining components $\Omega_{n+1},\Omega_{n+2}$, ``inside and outside" the chain of circles. We write $\Omega = \Omega_1 \sqcup \ldots \sqcup \Omega_{n+2}$.
    See \Cref{fig:2complex}.
    Since $\gamma$ transversely intersects $\check \Sigma_e^\circ$, the endpoints at infinity of $\widetilde \gamma$ are in distinct components of $\Omega$.

 \begin{figure}[h]
    \centering
    \begin{tikzpicture}
        \node at (0, 0) {\includegraphics[width=\linewidth]{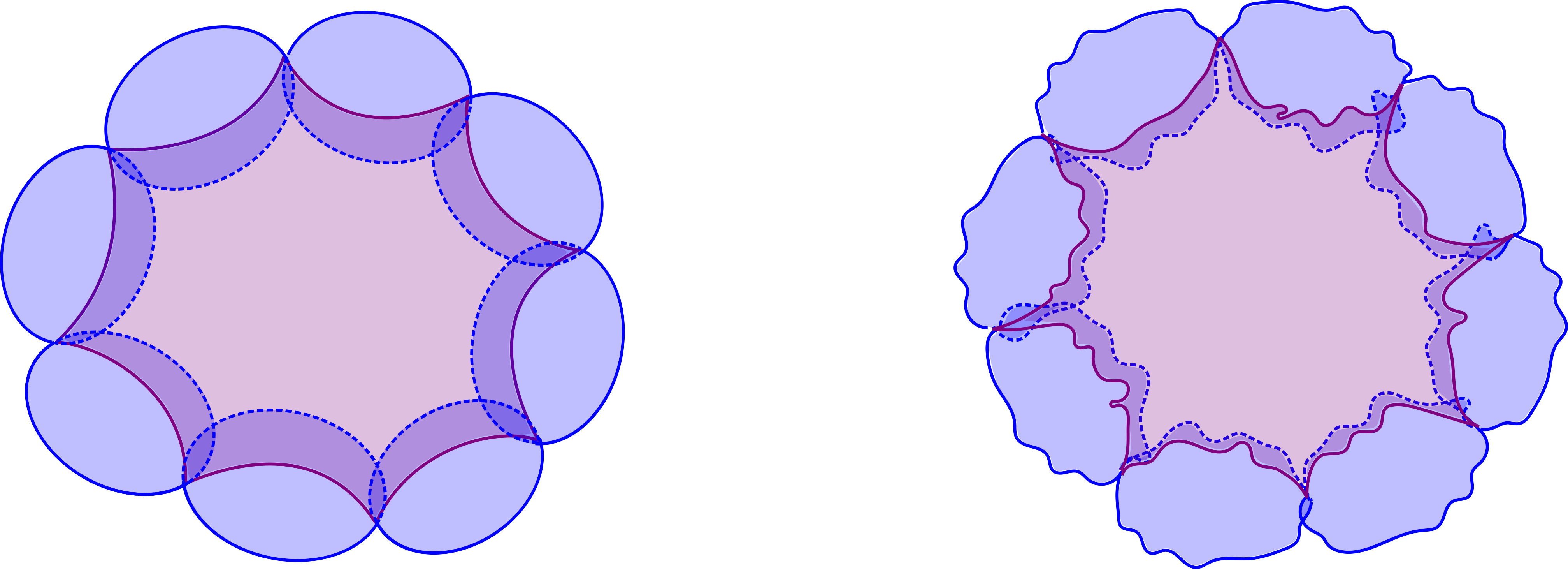}};
        \node at (-4, 0) {\color{violet}$\tilde{\Sigma}^{\circ}_e$};
        \node at (4.25, 0) {\color{violet}$\tilde{\Sigma}^{h}_e$};
        \node at (0, 0) {$\longrightarrow$};
        \node at (-4, -2.75) {$\widetilde{M}_{L'}$};
        \node at (4.25, -2.75) {$\widetilde{M}^h_{L'}$};
        \node at (-6, 2) {\color{blue}$\partial\Omega_i$};
        \node at (1.9, 2) {\color{blue}$\partial\Omega^h_i$};
    \end{tikzpicture}
    \caption{Illustration of the chain of circles in $\widetilde{M}_{L'}$ and quasi-circles in $\widetilde{M}^h_{L'}$.}
    \label{fig:2complex}
\end{figure}   

    Again by \Cref{T:CAT-1 to hyp}, the identity lifts to a quasi-isometry $\widetilde M_{L'} \to \widetilde M_{L'}^h \cong \mathbb H^3$.  There is an elevation $\widetilde \gamma^h$ of $\gamma^h$ to $\widetilde M_{L'}^h$ which is a bounded Hausdorff distance to the image of $\widetilde \gamma$.  Similarly, each of the surfaces in $\widetilde M_{L'}$ above is sent within a bounded Hausdorff distance of a pleated surface covering $\check \Sigma_e^h \to M_{L'}^h$ or $\Sigma_{f\pm}^h \to M_{L'}^h$.  The latter pleated surfaces are pleated hyperbolic planes with (quasi-)circles at infinity in $\partial \widetilde M_{L'}^h$ that are images of the circles in $\partial \widetilde M_{L'}$.  The complement of the union of these circles is the image of $\Omega$, and hence has $n+2$ components.  Furthermore, since the endpoints of $\widetilde \gamma^h$ are the same as those of the image of $\widetilde \gamma$, these endpoints are in different components.
    
    Next, observe that pleated surface covering $\check \Sigma_e^h$ is a pleated ideal hyperbolic $n$--gon with ideal endpoints at the intersection points of the circles.  The sides of the $n$--gon are geodesics in the pleated surfaces that cover the face surfaces.  In particular, the union of these $n+1$ pleated surfaces separate the $n+2$ regions at infinity in $\widetilde M_{L'}^h$, and thus $\widetilde \gamma^h$ must intersect at least one of them.  If it intersects one of the surfaces covering a face surface, then $\gamma^h$ intersects that face surface, and so comes within distance $0$ of it, proving the lemma in this case.  Otherwise,  $\gamma^h$ intersects the pleated cut edge surface.  Since this is a pleated ideal $n$--gon with $n \leq 20$, its area is at most $18\pi$.  Therefore, the minimal distance $\rho$ from a point of intersection of $\gamma^h$ with $\check \Sigma_e^h$ to a side of that polygon satisfies
    \[ \mbox{Area}(B_h(\rho)) = 4 \pi \sinh^2(\rho/2) \leq 18 \pi,\]
    where $B_h(\rho)$ is a hyperbolic ball of radius $\rho$.  Consequently, $\rho \leq 2\sinh^{-1}(3/\sqrt{2}) < 3$. 
    Therefore, $\gamma^h$ comes within distance $3$ of some $\Sigma_{f\pm}^h$ for $f \in \F(e)$, as required. 
\end{proof}

\begin{proof}[Proof of Proposition \ref{prop:predehnfillinglongessentialsys}]
Recall that we have fixed $R_0 > 0$, and want to show that $\esssys(L') \geq R_0$, provided $\combsys(\tau)$ is sufficiently large, where $L' = L'(\tau)$.  How large we must take $\combsys(\tau)$ will depend on other constants, which we now choose.  Without loss of generality, we assume $R_0 > 3$.

Let $\delta_0 > 0$ denote any upper bound on the diameter of the $\epsilon_0$--thick part of a hyperbolic $4$--punctured sphere (which exists by an area argument).
Next, let $\epsilon < \epsilon_0$ be any positive constant so that $R(\epsilon) > R_0 + \epsilon_0 + \delta_0$.  Similar to Lemma~\ref{lem:thinpart}, the resuls of \cite{bm,meyerhoff} imply that there is some positive $\epsilon' < \epsilon$, so a point in any hyperbolic $3$--manifold which is within distance $3$ of the $\epsilon$--thick part is in the $\epsilon'$--thick part.  We also let $\delta > 0$ denote any upper bound on the diameter of the $\epsilon'$--thick part of a $4$--punctured sphere.

Suppose $K = R_0 + 3+\delta$ and let $\eta_0 \geq 2$ be a natural number as in Lemma~\ref{lem:bounded abelian} for this value of $K$ and $\epsilon$.   Suppose $\tau$ is a Cooper-Thurston triangulation with 
\[ \left\lfloor \frac{\combsys(\tau)}{16} \right\rfloor = \eta \geq \eta_0,\]
and let $L' = L'(\tau)$ be the associated link.  We suppose $\esssys(L') \leq R_0$, and derive a contradiction.  This will prove the proposition.

Let $\gamma^h$ be a hyperbolic geodesic in $M_{L'}^h$ that realizes $\esssys(L')$, and $\gamma$ in $M_{L'}$ the homotopic geodesic in the locally $\CAT(-1)$ metric.  Let $T_1,\ldots,T_\eta$ be the tetrahedra from Lemma~\ref{lem:many tetrahedra} which $\gamma$ intersects.  By Lemma~\ref{lem:tet intersect}, $\gamma$ has non-empty transverse intersection with some surfaces in $\Sigma_\tau^\circ(T_i)$, for all $i=1,\ldots,\eta$.  By Lemmas~\ref{lem:meet the tet pleat} and \ref{lem:near a surface}, there exist pleated surfaces $\Sigma_{x_1}^h,\ldots,\Sigma_{x_\eta}^h \to M_{L'}^h$ which comes within distance $3$ of a point of $\gamma^h$, where one of the following holds:
\begin{enumerate}
    \item $x_i$ is a tetrahedron $T_i'$ where $T_i \cap T'_i \neq \varnothing$ or 
    \item $x_i = f_i\pm$, and $f_i$ is a face of some $T_i'$, where $T_i'$ and $T_i$ share an edge $e_i$. 
\end{enumerate}
In particular, in both cases, we have $d_{\tau}(T_i, T_i') \leq 2$. 
Since $T_1,\ldots,T_\eta$ are as in Lemma~\ref{lem:many tetrahedra}, the triangle inequality implies $d_\tau(T_i',T_j') \geq 4$, for all $i \neq j$.  By considering the totally geodesic representatives of the given surfaces in $M_{L'}$, we see that the only loops in any one of these surfaces that can be freely homotoped into another are the peripheral loops in the cusps of face surfaces, for cusps corresponding to the punctures coming from the focal component $K' \subset L'$; see Definition~\ref{def:main link}.  We call these the {\em focal cusps} of the pleated face surfaces.
We let $z_i \in \Sigma_{x_i}^h$ and $y_i \in \gamma^h$ be points with $d(y_i,z_i) \leq 3$, for all $i = 1, 2,\ldots, \eta$.

We divide the proof into two cases, depending on whether or not $\gamma^h$ meets the $\epsilon$--thin part of $M_{L'}^h$.  First, suppose $\gamma^h$ meets the $\epsilon$--thin part.  Since we assumed that $\ell(\gamma^h) \leq R_0$, Lemma~\ref{lem:where is esssys} implies $\gamma^h$ is contained in the $\epsilon$--thin part, and since $\gamma^h$ realizes the essential systole, it is the core geodesic of a $\epsilon$--thin Margulis tube, $\mathbb T_\epsilon$.
Since $d(y_i,z_i) \leq 3 < R(\epsilon)$,  Lemma~\ref{lem:thinpart} implies that $z_i \in \Sigma_{x_i}^h$ is in the component $\mathbb T_{\epsilon_0}$ of the $\epsilon_0$--thin part of $M_{L'}^h$ with $\mathbb T_\epsilon \subset T_{\epsilon_0}$.  Within distance $\delta_0$ of $z_i$ is a point in some $\epsilon_0$--thin part, hence some loop $\alpha_i \subset \Sigma_{x_i}^h$, of length at most $2 \epsilon_0$.  Every point on $\alpha_i$ is within distance $3+\delta_0 + \epsilon_0 < R(\epsilon)$ of $y_i$.  Consequently, $\alpha_i$ is contained in $T_{\epsilon_0}$, and is thus homotopic to a power $\gamma^h$ (the core of $\mathbb T_{\epsilon_0}$), for all $i = 1,\ldots, \eta$.  This is a contradiction since the only nontrivial loops in distinct pleated surfaces in the set $\Sigma_{x_1}^h,\ldots,\Sigma_{x_\eta}^h$ that can be homotopic to each other are loops in focal cusps.
Therefore, we can assume that $\gamma^h$ is completely contained in the $\epsilon$--thick part of $M_{L'}^h$.

Since $z_i$ is distance at most $3$ from $y_i$, it follows from our choice of constants that $z_i$ is in the $\epsilon'$--thick part of $M_{L'}^h$, and hence the $\epsilon'$--thick part of $\Sigma_{x_i}^h$.  Consequently, $z_i$ has distance at most $\delta$ from some point $w_i$ in a component of the $\epsilon'$--thin part of $\Sigma_{x_i}^h$, which we assume (as we may) is not a focal cusp.  In particular, $w_i$ is contained in a component, $\mathbb T_{\epsilon}^i$, of the $\epsilon$--thin part of $M_{L'}^h$.  Since these are not focal cusp thin parts, $\mathbb T_{\epsilon}^1,\ldots,\mathbb T_{\epsilon}^\eta$ are all distinct thin parts.  The distance $y_i$ to $w_i$ is at most $3 + \delta$, and so picking some point $y_0 \in \gamma^h$, we see that $d(y_0,w_i) \leq R_0+3+\delta = K$.  Choosing paths from $y_0$ to $w_i$ whose length is at most $K$, we can lift these to paths from a point $\widetilde y_0 \in \mathbb H^3$ to points $\widetilde w_1,\ldots,\widetilde w_\eta \in \mathbb H^3$.  The points $\widetilde w_1,\ldots,\widetilde w_\eta$ all lie within distance $K$ of $\widetilde y_0$ and are in distinct $\epsilon$--thin parts for $\Gamma = \pi_1(M_{L'}^h)$.  By Lemma~\ref{lem:bounded abelian}, $\eta < \eta_0$, which is a contradiction.
Since both cases are impossible, it must be that $\gamma^h$ has length greater than $R_0$, as required.
\end{proof}

\section{Filling knots}\label{sec:fillingknots}

In this final section, we prove Theorem~\ref{thm:main2specific} about the knots $K(n,\tau)$ (see \S\ref{sec:constructingknots}) which in particular proves the existence of filling knots (Theorem~\ref{T:filling knot}).  Like in \cite{lw}, we apply \Cref{thm:KW}, though we cannot directly appeal to Theorem~\ref{thm:main1specific}, since the theorem guarantees large essential {\em geodesic} systole, rather than simply large essential systole.  In particular, the action of $\pi_1(M \sss K(n,\tau))$ on hyperbolic $3$--space may be parabolic on elements that are non-trivial in $\pi_1(M)$.  The action of $\pi_1(M \sss K(n,\tau))$ on a hyperbolic space for which we apply \Cref{thm:KW} will be built from the locally $\CAT(-1)$ metric on $M \sss L'(\tau)$ and a construction of Gromov and Thurston (see \cite{BleilerHodgson}).  We begin by briefly recalling this construction.

\subsection{Negatively curved Dehn filling} An important ingredient in our construction is the \emph{Gromov-Thurston $2\pi$--Theorem}, or more precisely, it's proof, which we now briefly recall; see \cite{BleilerHodgson}.  Given a horoball cusp $U \cong T^2 \times [0,\infty)$, we equip the boundary $\partial U$ with its natural Euclidean metric.  If $\alpha$ in $\partial U \cong T^2$ is a slope, we define its {\em length} to be the length of its geodesic representative.  The {\em transverse distance of $\alpha$} is the distance between the two boundary components of the Euclidean annulus $\partial U \sss \alpha$.

If $\alpha$ is a slope on $\partial U$ of length greater than $2\pi$, then there is a negatively curved metric on a solid torus $V \cong S^1 \times D^2$ which is isometric to $U$ in a neighborhood of the boundary, with the property that the slope of the meridian $\{\ast\} \times \partial D^2 \subset \partial V$ is mapped by this isometry to $\alpha$ in $\partial U$.  Moreover, the metric is {\em maximally symmetric}: the ``obvious" isometric action of $O(2) \times O(2)$ on $\partial V \cong S^1 \times \partial D^2$ extends to an isometric action on $V$.  The fixed point sets of the first-factor-reflections define a foliation by totally geodesic meridian disks of $V$.  
The construction is quite flexible, and a minor modification of the proof of \cite[Lemma~4]{Leininger06} implies the following.


\begin{lemma} \label{lemma:technical filling}
    There exists $r > 2 \pi$ with the following property.  If $U$ is a standard horoball cusp and $\alpha$ is a slope on $\partial U$ with length 
    $\ell > r$, then there is a maximally symmetric Riemannian metric on a solid torus $V$ with sectional curvature in $[-2,-1/2]$, with the following properties:
    \begin{enumerate}
    \item There is an isometry $h \colon {\mathcal N}_1(\partial U) \to {\mathcal N}_1(\partial V)$ of the $1$--neighborhoods of the boundaries such that $h(\alpha)$ is a meridian in $\partial V$, and
    \item there is a constant $\xi$, depending only on $\ell$, so that the length of the core geodesic of $V$ is $\xi$ times the transverse distance of $\alpha$. \qed
    \end{enumerate} 
\end{lemma}

We will use Lemma~\ref{lemma:technical filling} to carry out a \emph{Gromov--Thurston negatively curved Dehn surgery} on a manifold $N$ with a locally $\CAT(-\tfrac12)$, piecewise Riemannian (or hyperbolic) metric containing a hyperbolic cusp $U \subset N$. To do so, suppose we have a slope $\alpha \subset \partial U$ of length $\ell > r$, where $r$ is the constant from Lemma~\ref{lemma:technical filling}. Take a solid torus $V$ given by the lemma.  Now we remove $U \sss {\mathcal N}_1(\partial U)$ from $N$, and glue in $V$ to the result via the isometry $h \colon {\mathcal N}_1(\partial U) \to {\mathcal N}_1(\partial V)$.  The resulting piecewise Riemannian metric on the Dehn surgered manifold $N(\alpha)$ is again locally $\CAT(-\tfrac12)$ and agrees outside of $V$ with the original piecewise Riemannian metric on $N$ outside of $U$; that is, $N(\alpha) \sss V$ is isometric to $N \sss U$ (by construction).

We will also need an orbifold version of the lemma and this Dehn surgery construction.  For this, we consider a cusp $U$ as before and a positive integer $q > 1$.  Now suppose $\alpha$ is a slope on $\partial U$ with length $\ell$ such that $q\cdot \ell > r$.  Consider a $q$--fold cyclic cover of $\widetilde U \to U$ corresponding to a homomorphism $\pi_1(U) \to \mathbb Z/q \mathbb Z$ for which $\alpha$ is sent to a generator.  Note that $\alpha^q$ lifts to a slope $\widetilde \alpha$ in $\partial \widetilde U$ with length $q \cdot \ell > r$.  Let $\widetilde V$ be the solid torus from the lemma applied to $\widetilde U$ and $\widetilde \alpha$.  Also let $\widetilde h \colon {\mathcal N}_1(\partial \widetilde U) \to {\mathcal N}_1(\partial \widetilde V)$ be the isometry.

\begin{figure}[H]
    \centering
    \begin{tikzpicture}
    \node at (0,0) {\includegraphics[width=\linewidth]{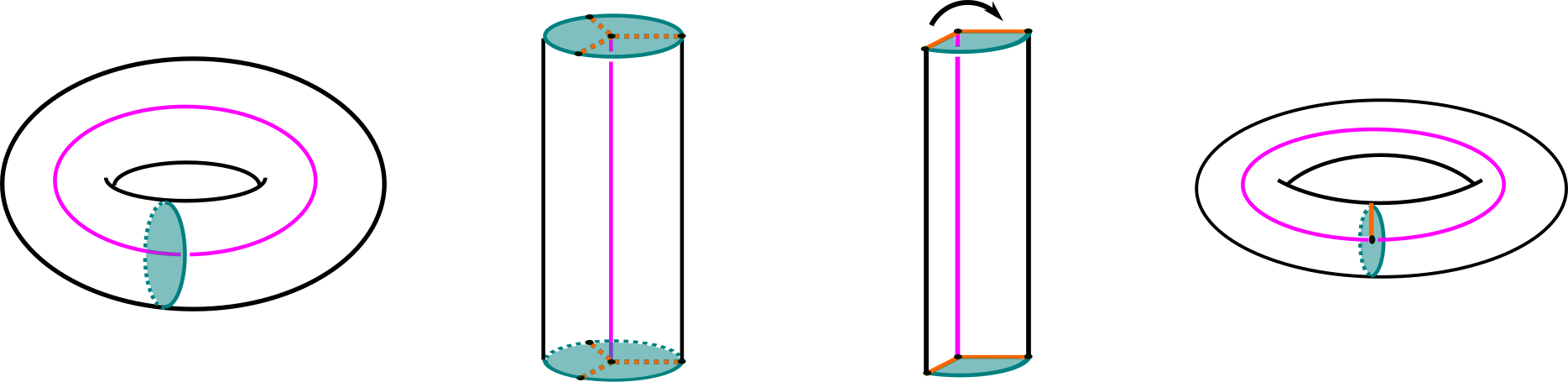}};
    \node at (-5, -1.5) {$U$};
    \node at (5, -1.5) {$V$};
    \node at (-2.7, 0.5) {cut};
    \node at (-2.7, 0) {$\longrightarrow$};
    \node at (0, 0.5) {$\mathbb{Z}/q\mathbb{Z}$};
    \node at (0, 0) {$\longrightarrow$}; 
    \node at (2.7, 0.5) {glue};
    \node at (2.7, 0) {$\longrightarrow$}; 
    \node at (1.75,0) {\color{mypink}{$3$}};
    \end{tikzpicture}
    \caption{Metric near the orbifold locus. Here, $q = 3$, so the core curve of $V$ is a singular curve of order $3$.}
    \label{fig:orbifoldlocus}
\end{figure}

The action of $\mathbb Z/q \mathbb Z$ on $\mathcal U$ with quotient $\widetilde U \to U$ determines an action on $\widetilde V$ via conjugation by $h$ (and extending): indeed, $\widetilde \alpha$ is invariant by the action, and $h(\widetilde \alpha)$ bounds a meridian disk in $\widetilde V$.  Let $\widetilde V \to V$ denote the quotient by this action.  The quotient $V$ is also a solid torus {\em orbifold} with Riemannian  metric with curvature in $[-2,-\frac12]$ and orbifold locus the geodesic core curve for $V$, with local group $\mathbb Z/q \mathbb Z$ (see e.g.~\cite[Chapter 13]{Thurstonnotes}).  The map $\widetilde h$ descends to an isometry $h \colon {\mathcal N}_1(\partial U) \to {\mathcal N}_1(\partial V)$. See \Cref{fig:orbifoldlocus}. 

As with the Gromov--Thurston negatively curved Dehn surgery above, given a locally $\CAT(-\tfrac12)$, piecewise Riemannian manifold $N$ containing a hyperbolic cusp $U$, we can perform Gromov--Thurston negatively curved $q$--orbifold surgery, for a slope $\alpha$ of length $\ell$ on $\partial U$, as long as $q \cdot \ell > r$.  The result in this case is a locally $\CAT(-\tfrac12)$, piecewise Riemannian orbifold, whose underlying space is $N(\alpha)$ and which has orbifold locus a geodesic core of the filling solid torus, with local group $\mathbb Z/q \mathbb Z$.

\subsection{Constructing an appropriate action}

Recall from Lemma~\ref{lem:final new link} that, given any Cooper--Thurston triangulation $\tau$, we have an associated link $L' = L'(\tau) \subset M$.  The complement $M_{L'}$ is equipped with a complete, piecewise hyperbolic, locally $\CAT(-1)$ metric. Furthermore, from \eqref{eq:L and L' complement}, $M_{L'}$ is obtained by isometrically gluing copies of the model tangle complement $T_0 \sss J_0$ (see \S\ref{subsection:fillinglinkrecap}) according to $\tau$, with the connector $\mathbb B^f$ inserted between tetrahedra glued along the face $f$, for each $f \in \F_\tau$.

We also recall from Lemma~\ref{lem:final new link} that the link $L'$ contains a component $K' = K'(\tau) \subset L'$ that meets every tetrahedron $T \in \T_\tau$ in the union of the $6$ arcs, and each connector $\mathbb B^f$ in a union of $3$ arcs.  We choose a fixed (pairwise disjoint) union of cusp neighborhoods $U_0 \subset T_0 \sss J_0$ of the ends corresponding to these arc components, and let $U_T \subset T \sss L'$ be the image under the isometry $T \sss L' \cong T_0 \sss J_0$. For each $f \in F_{\tau}$, we similarly choose pairwise disjoint cusp neighborhoods of the ends of each connector $U^f \subset \mathbb B^f$, for the ends corresponding to the three arcs in $\mathbb B^f$.  
There are only 8 isometry types of connectors (one for each choice of whether the twist in the disk $D_{fij}$ is performed or not in the construction from Lemma~\ref{lem:final new link}).  We choose the cusps $U^f \subset \mathbb B^f$ to be such that the isometry to one of the ``8 model connectors" sends $U^f$ to some fixed cusp of the model.  Finally, we make all these choices (in the models, independent of $M$ and $\tau$) so that upon gluing, we obtain a cusp neighborhood of the end of $M_{L'}$ corresponding to $K'$.  We denote this cusp neighborhood $U_{K'} \subset M_{L'}$.  For any $R > 0$, we let $U_{K'}^R \subset U_{K'}$ be the smaller cusp whose distance to $\partial U_{K'}$ is $R$.  Let $\mu$ denote the meridian on $\partial U_{K'}$, i.e.~the isotopy class that bounds a disk in $M$.  We will view $\mu$ as a slope on both $\partial U_{K'}$ and $\partial U_{K'}^R$, for any $R$.  Observe that $\mu$ is the the isotopy class of any component of intersection of $\partial U_{K'}$ with any face surface $\Sigma_{f\pm}^\circ$.  Since $\Sigma_{f\pm}^\circ$ is totally geodesic, its intersection with $\partial U_{K'}$ consists of geodesic representatives of the isotopy class in the flat metric, and its length is the length of $\mu$ in $\partial U_{K'}$.  In particular, this length, which we denote by $\ell_0$, is independent of $M$ and $\tau$.  The length of $\mu$ in $\partial U_{K'}^R$ is also independent of $M$ and $\tau$, and in fact is just given by $e^{-R}\ell_0$, since $\partial U_{K'}^R$ is isometric to $\partial U_{K'}$, scaled by $e^{-R}$.

\begin{lemma} \label{lem:uniform transverse distance}
    There exists a constant $C > 0$ with the following property.  For any $M$ and $\tau$, the transverse distance of $\mu$ in $\partial U_{K'}^R$ is at least
    \[ e^{-R}C|\T_\tau|.\]
    Moreover, the length of $\alpha$ on $U_{K'}^R$ is $e^{-R}\ell_0$.
\end{lemma}
\begin{proof}
    We already observe the last statement before the statement of the lemma, so we need only prove the first one.
    
    Let $\partial U_0$ denote the horospherical part of the boundary of $U_0$ in $T_0 \sss J_0$ (i.e.~it does not include the part of $\partial U_0$ contained in $\partial T_0 \sss J_0$). We observe that $\partial U_0 \cap (T_0 \sss J_0)$ is the disjoint union of $6$ isometric Euclidean annuli.  Let $c > 0$ denote the distance between the boundary components inside any one of these annuli.    
    Since $U_{K'}$ meets $T \sss L'$ in a cusp isometric to $U_0$, it follows that every tetrahedron contributes $6c$ to the transverse distance of $\mu$ in $\partial U_{K'}$.  Setting $C = 6c$, and again noting that $\partial U_{K'}^R$ is isometric to $\partial U_{K'}$ scaled by $e^{-R}$, completes the proof.
\end{proof}

\begin{lemma} \label{lem:separated tetrahedra}
    Given $R \geq 0$, there exists $\varepsilon > 0$ with the following property.  Suppose $M$ is any closed $3$--manifold and $\tau$ any Cooper--Thurston triangulation.  Then for any two tetrahedra in $T,T' \in \T_\tau$, with $d_\tau(T,T') \geq 4$, we have
    \[d_{M_{L'}\sss U_{K'}^R}(x,y) > \varepsilon\]
    for all $x \in T \sss (L' \cup U_{K'}^R)$ and $y \in T' \sss (L' \cup U_{K'}^R)$.
\end{lemma}
That is, in $M_{L'} \sss U_{K'}^R$, the distance between points in disjoint tetrahedra is at least $\varepsilon$.
\begin{proof}
    Except for the end corresponding to $K'$, all other ends have cusp neighborhoods either contained in a single connector, or meeting a single connector and the two tetrahedra to which it is glued (note that these tetrahedra are adjacent).  Thus, after removing $U_{K'}$, for any tetrahedron $T \in \T_\tau$, there is a $\varepsilon_T > 0$ so that the $\varepsilon_T$--neighborhood of $T \sss (L' \cup U_{K'}^R)$ in $M_{L'} \sss U_{K'}^R$ only meets adjacent tetrahedra and adjacent connectors.  Since the link of any vertex of a Cooper--Thurston triangulations is isomorphic to one of five $2$--complexes, and since there are only $8$ isometry types of connectors, it follows that the isometry type of
    \[ \left( \bigcup_{d_\tau(T,T') \leq 2} T' \right) \sss (L' \cup U_{K'}^R)
    \]
    is one of finitely many types, independent of $M$ and $\tau$.  In particular, we can choose $\varepsilon_T$ depending only on the isometry type the above submanifold.  Doing so, there is a uniform lower bound $\varepsilon' > 0$ to such $\varepsilon_T$, and setting $\varepsilon = \varepsilon'/2$, proves the lemma.
\end{proof}

A non-null-homotopic loop in $M_{L'}$ must pass through $\lfloor \frac{\combsys(\tau)}{16}\rfloor$ tetrahedra that are pairwise distance at least $8$ apart, by Lemma \ref{lem:many tetrahedra}.  In particular, these tetrahedra are disjoint, and so such a curve which is contained in $M_{L'} \sss U_{K'}^R$ must travel distance $\varepsilon$ in passing from one tetrahedron to the next.  Thus, combined with Lemma~\ref{lem:separated tetrahedra} we immediately find the following. 

\begin{corollary}
    If $\gamma$ is a closed curve in $M_{L'} \sss U_{K'}^R$, which is non-null-homotopic in $M$, then its length in the locally $\CAT(-1)$ metric is bounded below by $\varepsilon \cdot \lfloor \frac{\combsys(\tau)}{16}\rfloor$. \qed
\end{corollary}

Given $M$, $\tau$, and $n \in \mathbb Z$, we write $M_{K(n,\tau)} = M \sss K(n,\tau)$ (see \S\ref{sec:constructingknots}).

\begin{corollary} \label{cor:long curves cat-1/2}
    For any Cooper--Thurston triangulation $\tau$ of a closed $3$--manifold $M$, and $R \geq 0$, there exists $n_0$ so that for all $n$ with $|n| > n_0$, $M_{K(n,\tau)}$ admits a $\CAT(-\tfrac12)$ metric with the following properties.
    \begin{enumerate}
        \item There is a cusp neighborhood $U_n^R$ in $M_{K(n,\tau)}$ and isometry $U_n^R \cong U_{K'}^R \subset M_{L'}$ for which the meridian $\mu_n$ is the image of $\mu$ under the isometry, and
        \item if $\gamma$ is a closed curve in $M_{K(n,\tau)} \sss U_n^R$ which is non-null-homotopic in $M$, then the length of $\gamma$ is greater than $\varepsilon \cdot \lfloor \tfrac{\combsys(\tau)}{16} \rfloor$.
    \end{enumerate}
\end{corollary}
Note that the case $R = 0$ gives us $U_n^0$ is isometric to $U_{K'}^0 = U_{K'}$.  We write $U_n = U_n^0$.
\begin{proof}
    Choose pairwise disjoint cusps $W_1,\ldots,W_k$ for the ends of $M_{L'}$, different than the one corresponding to $K'$, which are also disjoint from $U_{K'}$.  Inside of each, we choose smaller cusps $W_i' \subset W_i$, for all $i$, whose distance to $\partial W_i$ is at least $\varepsilon \cdot \lfloor \combsys(\tau)/16 \rfloor$.
    Choose $n_0 > 0$ so that for all $n$ with $|n| > n_0$ the slope $(1,n)$ on $\partial W_i'$ has length at least $r$, where $r$ is the constant from Lemma~\ref{lemma:technical filling}.
    Now perform the $(1,n)$ negatively curved Dehn filling on each $W_i'$ (as explained in the discussion following the Lemma~\ref{lemma:technical filling}).  Suppose $V_i$ is the filling solid torus from the lemma.  The result is a locally $\CAT(-\tfrac12)$ metric on $M_{K(n,\tau)}$, which outside the union of $V_i$ is isometric to $M_{L'}$, outside the union of the $W_i'$.  In particular, since all the cusps $U_{K'}^R$ are disjoint from all $V_i$, the given isometry restricts to an isometry from a cusp neighborhood $U_n^R \subset M_{K(n,\tau)}$ to $U_{K'}^R$.  Since $K(n,\tau)$ is the image of $K'$ by the homeomorphism $M \to M$ as described in \S\ref{sec:constructingknots}, the meridian $\mu_n$ in $\partial U_n^R$ and $\mu$ are related by the isometry, proving the first claim.

    Now suppose $\gamma$ in $M_{K(n,\tau)} \sss U_n^R$ is any curve which is non-null-homotopic in $M$. We observe that if $\gamma$ meets any of the solid tori $V_i$, then since it cannot be contained in any larger solid torus containing $V_i$ (since these are null homotopic in $M$), it must intersect $\partial W_i$ to get to $V_i$.  This requires distance at least $\epsilon \cdot \lfloor \tfrac{\combsys(\tau)}{16} \rfloor$ by our choice of cusps $W_i$, and so $\gamma$ must have length at least this number, as required.  Otherwise, $\gamma$ is disjoint from the filling solid tori, hence it is related by an isometry to a closed curve in $M_{L'} \sss U_{K'}^R$.  By Corollary~\ref{cor:long curves cat-1/2}, the length must be at least $\epsilon \cdot \lfloor \frac{\combsys(\tau)}{16} \rfloor$, completing the proof of the second part, hence the proof of the corollary.
\end{proof}

We now describe the key proposition needed to prove the theorem.

\begin{proposition} \label{prop:final}
Suppose $M$ is closed orientable $3$--manifold.  Given $R > 0$, there exists $\tau$ and $n_0$, such that for all $n$ with $|n|>n_0$, there is an action of $\pi_1(M_{K(n,\tau)})$ on a $\CAT(-\tfrac12)$ space, $X$.  Under this action, every element of $\pi_1(M_{K(n,\tau)})$ which is non-null-homotopic in $M$ acts loxodromically with translation length at least $R$.
\end{proposition}

\begin{proof}[Proof of Proposition~\ref{prop:final}]
Given $R > 0$, let $\varepsilon > 0$ be as in Lemma~\ref{lem:separated tetrahedra}.  Let $C$ and $\ell_0$ be the constants from Lemma~\ref{lem:uniform transverse distance}.  Fix $q > 0$ large enough so that $q \cdot e^{-R} \ell_0 > r$, where $r > 0$ is the constant from Lemma~\ref{lemma:technical filling}.  Let $\xi$ be the constant from Lemma~\ref{lemma:technical filling}, depending on the length, $q \cdot e^{-R} \ell_0 > r$.

Suppose $\tau$ is any Cooper--Thurston triangulation of $M$ with
\[ \min\left\{\xi e^{-R}C|\T_\tau|,\epsilon \cdot \lfloor \tfrac{\combsys(\tau)}{16} \rfloor\right\}  > R \] and suppose $n_0 = n_0(R,\tau)$ from Corollary~\ref{cor:long curves cat-1/2}.  Let $n$ be any integer with $|n| > n_0$ and $U_n^R \subset M_{K(n,\tau)}$ the cusp neighborhood from that corollary.
Since $U_n^R$ is isometric $U_{K'}^R$ (also from that corollary), Lemma~\ref{lem:uniform transverse distance}, implies that the transverse distance of $\mu_n$ is $e^{-R}C|\T_\tau|$ and the length of $\mu_n$ is $e^{-R}\ell_0$.
Let $V$ denote the $q$--orbifold solid torus for which there is an isometry ${\mathcal N}_1(\partial V) \cong {\mathcal N}_1(\partial U_n^R)$ for the slope $\mu_n$ as in Lemma~\ref{lemma:technical filling}.  By the lemma, the length of the core geodesic, $K_q$, (i.e.~the orbifold locus) of $V$ given by
\[ \xi e^{-R}C|\T_\tau| > R. \]

Now, we perform the negatively curved $q$--orbifold Gromov--Thurston surgery along $\mu_n$ on the cusp $U_n^R$ of $M_{K(n,\tau)}$ with the solid torus $V$.  For simplicity, we denote the result $M_{K(n,\tau)}(\mu_n)_q$, simply as $N_q \cong M_{K(n,\tau)}(\mu_n)_q$ for the remainder of the proof.  Given a loop $\gamma$ in $M_{K(n,\tau)}$, we can take a sequence of representatives infimizing the lengths, and extract a limit, $\gamma^*$ in $N_q$ (which exists since $N_q$ is compact). We denote the length of $\gamma^*$ by $l(\gamma) = l(\gamma^*)$.

\begin{claim}
If $\gamma$ is any closed loop in $M_{K(n,\tau)}$ which is non-null-homotopic in $M$, then $l(\gamma) > R$.
\end{claim}
\begin{proof}
    Since $U_n^R \subset U_n$, it follows that $(U_n \sss U_n^R) \cup V$ is also an orbifold solid torus we denote $V' \subset N_q$.  If $\gamma^*$ is contained in $V'$, then $\gamma^*$ must be (a multiple of) the core geodesic $K_q$ of $V'$ (any other loop is null-homotopic in $M$).  This is also the core of $V$, and thus has length greater than $R$ as noted above. If $\gamma^*$ is not contained in $V'$, but intersects $V$, then it intersects both $\partial V'$ and $\partial V$, but the distance between these tori is $R$, hence $l(\gamma) > R$ again.  Therefore, we can assume that $\gamma^*$ is disjoint from $V$.  Then $\gamma^*$ is contained in $N_q \sss V$ which is isometric to $M_{K(n,\tau)} \sss U_n^R$.  By Corollary~\ref{cor:long curves cat-1/2}, $l(\gamma^*)$ is greater than $\epsilon \cdot \lfloor \tfrac{\combsys(\tau)}{16}\rfloor \geq R$, completing the proof of the claim.
\end{proof}
Now observe that $M_{K(n,\tau)}$ is homeomorphic to $N_q \sss K_q$.
Now consider the universal cover of $N_q \sss K_q$, equipped with the incomplete pull-back metric from its inclusion into $N_q$, and let $X$ denote the metric completion.  The action of $\pi_1(M_{K(n,\tau)}) = \pi_1(N_q \sss K_q)$ on the universal cover extends to $X$.  The quotient is precisely $N_q$.  Since $N_q$ is equipped with a locally $\CAT(-\tfrac12)$ orbifold metric, the universal cover of $N_q \sss K_q$ with the incomplete metric is locally $\CAT(-\tfrac12)$. Its metric completion, $X$, is obtained by adding bi--infinite geodesics, and is locally $\CAT(-\tfrac12)$ by \cite[Theorem~1.2]{allcock}. By the Cartan--Hadamard Theorem, $X$ is $\CAT(-\tfrac12)$ (see \cite[Theorem II.4.1]{bh}).

If $\gamma \in \pi_1(M_{K(n,\tau)})$ is non-null-homotopic in $M$, then $\gamma$ acts loxodromically: indeed, the axis projects to $\gamma^*$.  Consequently, the translation length is greater than $R$ by the claim, as required.
\end{proof}

With Proposition~\ref{prop:final} in hand, Theorem~\ref{thm:main2specific} now follows easily from \Cref{thm:KW}.

\begin{proof}[Proof of Theorem~\ref{thm:main2specific}]
Let $M$ be any closed $3$--manifold and $m > 0$, as in the statement of the theorem.  Let $C = \max\{C(k) \mid k \leq m \}$, where $C(k)$ is the constant from \Cref{thm:KW}.

By Proposition 1.2 of \cite[Section III.H.1]{bh}, any $\CAT(-\tfrac{1}{2})$ space is $\delta$--hyperbolic for some uniform $\delta > 0$.  Choose any $R > \delta C$.  By \Cref{prop:final}, there exists a Cooper--Thurston triangulation $\tau$ and $n \geq 0$, associated with a knot $K := K(n, \tau)$, such that $\pi_1(M_K)$ acts on a $\CAT(-\tfrac{1}{2})$ space $X$, and every loop which is non-null-homotopic in $M$ acts with translation length at least $R$ in $M_K$. 

Let $i: M_K \hookrightarrow M$ be the inclusion map. Consider a graph $\Gamma$ such that rank$(\pi_1(\Gamma)) = k \leq m$ and a map \[f: \Gamma \rightarrow M_K\] with the property that $i \circ f: \Gamma \rightarrow M$ has full rank (see \S\ref{sec:constructingknots}). Let $B = \langle h_1, \ldots, h_k\rangle$ be a basis for $\pi_1(\Gamma)$. Write $g_i := f_{\ast}(h_i)$ for all $1 \leq i \leq k$. Observe that $f_{\ast}(\pi_1(\Gamma))$ acts on $X$ by isometries, and since $i \circ f$ has full rank, $\langle g_1, \ldots, g_k\rangle$ is a generating set for $f_{\ast}(\pi_1(\Gamma))$. 

Suppose $\langle g_1, \ldots, g_k\rangle$ is Nielsen equivalent to some $\langle g'_1,\ldots, g'_k\rangle$, and \[d_{\tilde{X}}(g'_1.x, x) \leq \delta C \leq \delta C(k)\] for some $x \in X$. Then $i_{\ast}(g'_1)$ is trivial in $\pi_1(M)$ by our choice of $R$. By Lemma 2.1 of \cite{lw}, the Nielsen moves on $\langle g_1,\ldots,g_k\rangle$ lift to Nielsen moves on $\langle h_1, \ldots, h_k\rangle$, creating a new basis $\langle h'_1, \ldots, h'_k\rangle$ on $\pi_1(\Gamma)$ with the property that $f_{\ast}(h'_i) = g'_i$ for all $1 \leq i \leq k$. But then \[(i \circ f)_{\ast}(h'_1) = i_{\ast}(g'_1) = 0.\] This means $(i \circ f)(\pi_1(\Gamma))$, which has generating set $i_*(g_2'),\ldots,i_*(g_k')$, which has fewer than $k$ elements.  This contradicts the assumption that $i \circ f$ was full rank.

By Theorem~\ref{thm:KW}, $f_{\ast}(\pi_1(\Gamma))$ is free. By the Hopfian property for free groups, $f_{\ast}$ is injective. Since $f$ and $k \leq m$ were arbitrary, $K$ is $m$--rank filling in $M$. 
\end{proof}

\bibliographystyle{alpha}
\bibliography{main}

@article {adams,
    AUTHOR = {Adams, Colin C.},
     TITLE = {Thrice-punctured spheres in hyperbolic {$3$}-manifolds},
   JOURNAL = {Trans. Amer. Math. Soc.},
  FJOURNAL = {Transactions of the American Mathematical Society},
    VOLUME = {287},
      YEAR = {1985},
    NUMBER = {2},
     PAGES = {645--656},
      ISSN = {0002-9947,1088-6850},
   MRCLASS = {57N10 (57M25)},
  MRNUMBER = {768730},
       DOI = {10.2307/1999666},
       URL = {https://doi-org.ezproxy.rice.edu/10.2307/1999666},
}

@incollection {Adams-augment,
    AUTHOR = {Adams, Colin C.},
     TITLE = {Augmented alternating link complements are hyperbolic},
 BOOKTITLE = {Low-dimensional topology and {K}leinian groups
              ({C}oventry/{D}urham, 1984)},
    SERIES = {London Math. Soc. Lecture Note Ser.},
    VOLUME = {112},
     PAGES = {115--130},
 PUBLISHER = {Cambridge Univ. Press, Cambridge},
      YEAR = {1986},
      ISBN = {0-521-33905-7},
   MRCLASS = {57M25 (57N10)},
  MRNUMBER = {903861},
MRREVIEWER = {W.\ Heil},
}

@article {as,
    AUTHOR = {Adams, Colin and Schoenfeld, Eric},
     TITLE = {Totally geodesic {S}eifert surfaces in hyperbolic knot and
              link complements. {I}},
   JOURNAL = {Geom. Dedicata},
  FJOURNAL = {Geometriae Dedicata},
    VOLUME = {116},
      YEAR = {2005},
     PAGES = {237--247},
      ISSN = {0046-5755,1572-9168},
   MRCLASS = {57M25 (57M50)},
  MRNUMBER = {2195448},
MRREVIEWER = {Mattia\ Mecchia},
       DOI = {10.1007/s10711-005-9018-z}
}

@article {adamsreid,
    AUTHOR = {Adams, Colin C. and Reid, Alan W.},
     TITLE = {Systoles of hyperbolic {$3$}-manifolds},
   JOURNAL = {Math. Proc. Cambridge Philos. Soc.},
  FJOURNAL = {Mathematical Proceedings of the Cambridge Philosophical
              Society},
    VOLUME = {128},
      YEAR = {2000},
    NUMBER = {1},
     PAGES = {103--110},
      ISSN = {0305-0041,1469-8064},
   MRCLASS = {53C22 (53C20 57M50 57N10)},
  MRNUMBER = {1724432},
MRREVIEWER = {Athanase\ Papadopoulos},
       DOI = {10.1017/S0305004199003990},
       URL = {https://doi.org/10.1017/S0305004199003990},
}

@book {bh,
    AUTHOR = {Bridson, Martin R. and Haefliger, Andr\'e},
     TITLE = {Metric spaces of non-positive curvature},
    SERIES = {Grundlehren der mathematischen Wissenschaften [Fundamental
              Principles of Mathematical Sciences]},
    VOLUME = {319},
 PUBLISHER = {Springer-Verlag, Berlin},
      YEAR = {1999},
     PAGES = {xxii+643},
      ISBN = {3-540-64324-9},
   MRCLASS = {53C23 (20F65 53C70 57M07)},
  MRNUMBER = {1744486},
MRREVIEWER = {Athanase\ Papadopoulos},
       DOI = {10.1007/978-3-662-12494-9},
}

@article {bm,
    AUTHOR = {Brooks, Robert and Matelski, J. Peter},
     TITLE = {Collars in {K}leinian groups},
   JOURNAL = {Duke Math. J.},
  FJOURNAL = {Duke Mathematical Journal},
    VOLUME = {49},
      YEAR = {1982},
    NUMBER = {1},
     PAGES = {163--182},
      ISSN = {0012-7094,1547-7398},
   MRCLASS = {30F40 (53A35)},
  MRNUMBER = {650375},
MRREVIEWER = {Troels\ J\o rgensen},}

@article {FreeKrus2023,
    AUTHOR = {Freedman, Michael and Krushkal, Vyacheslav},
     TITLE = {Filling links and spines in 3-manifolds},
      NOTE = {With appendix by Christopher J. Leininger and Alan W. Reid},
   JOURNAL = {Comm. Anal. Geom.},
  FJOURNAL = {Communications in Analysis and Geometry},
    VOLUME = {31},
      YEAR = {2023},
    NUMBER = {10},
     PAGES = {2307--2333},
      ISSN = {1019-8385,1944-9992},
   MRCLASS = {57K30 (57K10)},
  MRNUMBER = {4785590},
MRREVIEWER = {Carlo\ Petronio},
       DOI = {10.4310/cag.2023.v31.n10.a1},
       URL = {https://doi.org/10.4310/cag.2023.v31.n10.a1},
}

@article {lakelandleininger,
    AUTHOR = {Lakeland, Grant S. and Leininger, Christopher J.},
     TITLE = {Systoles and {D}ehn surgery for hyperbolic 3-manifolds},
   JOURNAL = {Algebr. Geom. Topol.},
  FJOURNAL = {Algebraic \& Geometric Topology},
    VOLUME = {14},
      YEAR = {2014},
    NUMBER = {3},
     PAGES = {1441--1460},
      ISSN = {1472-2747,1472-2739},
   MRCLASS = {57M50 (57M27)},
  MRNUMBER = {3190600},
MRREVIEWER = {Yuichi\ Yamada},
       DOI = {10.2140/agt.2014.14.1441},
       URL = {https://doi.org/10.2140/agt.2014.14.1441},
}

@misc{lw,
      title={Filling Links and Essential Systole}, 
      author={Christopher J. Leininger and Yandi Wu},
      year={2025},
      eprint={2505.00896},
      archivePrefix={arXiv},
      primaryClass={math.GT},
      url={https://arxiv.org/abs/2505.00896}, 
}

@article {BleilerHodgson,
    AUTHOR = {Bleiler, Steven A. and Hodgson, Craig D.},
     TITLE = {Spherical space forms and {D}ehn filling},
   JOURNAL = {Topology},
  FJOURNAL = {Topology. An International Journal of Mathematics},
    VOLUME = {35},
      YEAR = {1996},
    NUMBER = {3},
     PAGES = {809--833},
      ISSN = {0040-9383},
   MRCLASS = {57M25 (57M50)},
  MRNUMBER = {1396779},
MRREVIEWER = {Ying\ Qing\ Wu},
       DOI = {10.1016/0040-9383(95)00040-2},
       URL = {https://doi.org/10.1016/0040-9383(95)00040-2},
}

@article {CoopThur1988,
    AUTHOR = {Cooper, D. and Thurston, W. P.},
     TITLE = {Triangulating {$3$}-manifolds using {$5$} vertex link types},
   JOURNAL = {Topology},
  FJOURNAL = {Topology. An International Journal of Mathematics},
    VOLUME = {27},
      YEAR = {1988},
    NUMBER = {1},
     PAGES = {23--25},
      ISSN = {0040-9383},
   MRCLASS = {57Q15},
  MRNUMBER = {935525},
MRREVIEWER = {Martin\ Scharlemann},
       DOI = {10.1016/0040-9383(88)90004-3},
       URL = {https://doi.org/10.1016/0040-9383(88)90004-3},
}

@article {Hoffman,
    AUTHOR = {Hoffman, Neil and Ichihara, Kazuhiro and Kashiwagi, Masahide
              and Masai, Hidetoshi and Oishi, Shin'ichi and Takayasu,
              Akitoshi},
     TITLE = {Verified computations for hyperbolic 3-manifolds},
   JOURNAL = {Exp. Math.},
  FJOURNAL = {Experimental Mathematics},
    VOLUME = {25},
      YEAR = {2016},
    NUMBER = {1},
     PAGES = {66--78},
      ISSN = {1058-6458,1944-950X},
   MRCLASS = {57N10 (57R05 65G40)},
  MRNUMBER = {3424833},
       DOI = {10.1080/10586458.2015.1029599},
       URL = {https://doi-org.ezproxy.rice.edu/10.1080/10586458.2015.1029599},
}

@article {kapovichweidmann,
    AUTHOR = {Kapovich, Ilya and Weidmann, Richard},
     TITLE = {Nielsen methods and groups acting on hyperbolic spaces},
   JOURNAL = {Geom. Dedicata},
  FJOURNAL = {Geometriae Dedicata},
    VOLUME = {98},
      YEAR = {2003},
     PAGES = {95--121},
      ISSN = {0046-5755,1572-9168},
   MRCLASS = {20F67 (20F65 57M07 57M60)},
  MRNUMBER = {1988426},
MRREVIEWER = {Eric\ M.\ Freden},
       DOI = {10.1023/A:1024064029186},
       URL = {https://doi.org/10.1023/A:1024064029186},
}

@article {Lackenby-augment,
    AUTHOR = {Lackenby, Marc},
     TITLE = {The volume of hyperbolic alternating link complements},
      NOTE = {With an appendix by Ian Agol and Dylan Thurston},
   JOURNAL = {Proc. London Math. Soc. (3)},
  FJOURNAL = {Proceedings of the London Mathematical Society. Third Series},
    VOLUME = {88},
      YEAR = {2004},
    NUMBER = {1},
     PAGES = {204--224},
      ISSN = {0024-6115,1460-244X},
   MRCLASS = {57M25 (57M50 57N10)},
  MRNUMBER = {2018964},
MRREVIEWER = {Carlo\ Petronio},
       DOI = {10.1112/S0024611503014291},
       URL = {https://doi.org/10.1112/S0024611503014291},
}

@article {Leininger06,
    AUTHOR = {Leininger, Christopher J.},
     TITLE = {Small curvature surfaces in hyperbolic 3-manifolds},
   JOURNAL = {J. Knot Theory Ramifications},
  FJOURNAL = {Journal of Knot Theory and its Ramifications},
    VOLUME = {15},
      YEAR = {2006},
    NUMBER = {3},
     PAGES = {379--411},
      ISSN = {0218-2165,1793-6527},
   MRCLASS = {57M50 (53C21 57M25)},
  MRNUMBER = {2217503},
MRREVIEWER = {Inkang\ Kim},
       DOI = {10.1142/S0218216506004531},
       URL = {https://doi.org/10.1142/S0218216506004531},
}

@article {meyerhoff,
    AUTHOR = {Meyerhoff, Robert},
     TITLE = {A lower bound for the volume of hyperbolic {$3$}-manifolds},
   JOURNAL = {Canad. J. Math.},
  FJOURNAL = {Canadian Journal of Mathematics. Journal Canadien de
              Math\'ematiques},
    VOLUME = {39},
      YEAR = {1987},
    NUMBER = {5},
     PAGES = {1038--1056},
      ISSN = {0008-414X,1496-4279},
   MRCLASS = {57S30 (22E40 53C20)},
  MRNUMBER = {918586},
MRREVIEWER = {Robert\ Brooks},
       DOI = {10.4153/CJM-1987-053-6},
       URL = {https://doi-org.ezproxy.rice.edu/10.4153/CJM-1987-053-6},
}

@article {nielsen,
    AUTHOR = {Nielsen, Jakob},
     TITLE = {Die {I}somorphismengruppe der freien {G}ruppen},
   JOURNAL = {Math. Ann.},
  FJOURNAL = {Mathematische Annalen},
    VOLUME = {91},
      YEAR = {1924},
    NUMBER = {3-4},
     PAGES = {169--209},
      ISSN = {0025-5831,1432-1807},
   MRCLASS = {99-04},
  MRNUMBER = {1512188},
       DOI = {10.1007/BF01556078},
       URL = {https://doi.org/10.1007/BF01556078},
}

@article {Purcell-augment,
    AUTHOR = {Purcell, Jessica S.},
     TITLE = {Hyperbolic geometry of multiply twisted knots},
   JOURNAL = {Comm. Anal. Geom.},
  FJOURNAL = {Communications in Analysis and Geometry},
    VOLUME = {18},
      YEAR = {2010},
    NUMBER = {1},
     PAGES = {101--120},
      ISSN = {1019-8385,1944-9992},
   MRCLASS = {57M50 (57M25)},
  MRNUMBER = {2660459},
MRREVIEWER = {Colin\ C.\ Adams},
       DOI = {10.4310/CAG.2010.v18.n1.a4},
       URL = {https://doi.org/10.4310/CAG.2010.v18.n1.a4},
}

@book {BenedettiPetronio,
    AUTHOR = {Benedetti, Riccardo and Petronio, Carlo},
     TITLE = {Lectures on hyperbolic geometry},
    SERIES = {Universitext},
 PUBLISHER = {Springer-Verlag, Berlin},
      YEAR = {1992},
     PAGES = {xiv+330},
      ISBN = {3-540-55534-X},
   MRCLASS = {57M50 (30F40 30F60 51M10 57N10)},
  MRNUMBER = {1219310},
MRREVIEWER = {Colin\ C.\ Adams},
       DOI = {10.1007/978-3-642-58158-8},
       URL = {https://doi-org.ezproxy.rice.edu/10.1007/978-3-642-58158-8},
}

@article {Thurston-hyp,
    AUTHOR = {Thurston, William P.},
     TITLE = {Hyperbolic structures on {$3$}-manifolds. {I}. {D}eformation
              of acylindrical manifolds},
   JOURNAL = {Ann. of Math. (2)},
  FJOURNAL = {Annals of Mathematics. Second Series},
    VOLUME = {124},
      YEAR = {1986},
    NUMBER = {2},
     PAGES = {203--246},
      ISSN = {0003-486X,1939-8980},
   MRCLASS = {57N10},
  MRNUMBER = {855294},
MRREVIEWER = {G.\ Peter\ Scott},
       DOI = {10.2307/1971277},
       URL = {https://doi.org/10.2307/1971277},
}

@article {McMullen,
    AUTHOR = {McMullen, Curt},
     TITLE = {Riemann surfaces and the geometrization of {$3$}-manifolds},
   JOURNAL = {Bull. Amer. Math. Soc. (N.S.)},
  FJOURNAL = {American Mathematical Society. Bulletin. New Series},
    VOLUME = {27},
      YEAR = {1992},
    NUMBER = {2},
     PAGES = {207--216},
      ISSN = {0273-0979,1088-9485},
   MRCLASS = {57M50 (30C75 30F60)},
  MRNUMBER = {1153266},
MRREVIEWER = {Colin\ C.\ Adams},
       DOI = {10.1090/S0273-0979-1992-00313-0},
       URL = {https://doi.org/10.1090/S0273-0979-1992-00313-0},
}

@incollection {Morgan,
    AUTHOR = {Morgan, John W.},
     TITLE = {On {T}hurston's uniformization theorem for three-dimensional
              manifolds},
 BOOKTITLE = {The {S}mith conjecture ({N}ew {Y}ork, 1979)},
    SERIES = {Pure Appl. Math.},
    VOLUME = {112},
     PAGES = {37--125},
 PUBLISHER = {Academic Press, Orlando, FL},
      YEAR = {1984},
      ISBN = {0-12-506980-4},
   MRCLASS = {57N10 (30F10)},
  MRNUMBER = {758464},
       DOI = {10.1016/S0079-8169(08)61637-2},
       URL = {https://doi.org/10.1016/S0079-8169(08)61637-2},
}

@article {allcock,
    AUTHOR = {Allcock, Daniel},
     TITLE = {Completions, branched covers, {A}rtin groups, and singularity
              theory},
   JOURNAL = {Duke Math. J.},
  FJOURNAL = {Duke Mathematical Journal},
    VOLUME = {162},
      YEAR = {2013},
    NUMBER = {14},
     PAGES = {2645--2689},
      ISSN = {0012-7094,1547-7398},
   MRCLASS = {53C23 (14B07 20F36 51K10 57N65)},
  MRNUMBER = {3127810},
MRREVIEWER = {Igor\ Belegradek},
       DOI = {10.1215/00127094-2380977},
       URL = {https://doi-org.ezproxy.rice.edu/10.1215/00127094-2380977},
}

@misc{Thurstonnotes,
  author        = {William P. Thurston},
  title         = {The geometry and topology of 3-manifolds},
  year          = {1978},
  publisher={Princeton University}
}

@article {thurston82,
    AUTHOR = {Thurston, William P.},
     TITLE = {Three-dimensional manifolds, {K}leinian groups and hyperbolic
              geometry},
   JOURNAL = {Bull. Amer. Math. Soc. (N.S.)},
  FJOURNAL = {American Mathematical Society. Bulletin. New Series},
    VOLUME = {6},
      YEAR = {1982},
    NUMBER = {3},
     PAGES = {357--381},
      ISSN = {0273-0979,1088-9485},
   MRCLASS = {57N10 (20H15 30F40 57M35 57M40 57S17)},
  MRNUMBER = {648524},
MRREVIEWER = {Klaus\ Johannson},
       DOI = {10.1090/S0273-0979-1982-15003-0},
       URL = {https://doi-org.ezproxy.rice.edu/10.1090/S0273-0979-1982-15003-0},
}

@article {white,
    AUTHOR = {White, Matthew E.},
     TITLE = {Injectivity radius and fundamental groups of hyperbolic
              3-manifolds},
   JOURNAL = {Comm. Anal. Geom.},
  FJOURNAL = {Communications in Analysis and Geometry},
    VOLUME = {10},
      YEAR = {2002},
    NUMBER = {2},
     PAGES = {377--395},
      ISSN = {1019-8385,1944-9992},
   MRCLASS = {57M50 (57N10)},
  MRNUMBER = {1900756},
MRREVIEWER = {Bruno\ P.\ Zimmermann},
       DOI = {10.4310/CAG.2002.v10.n2.a6},
       URL = {https://doi.org/10.4310/CAG.2002.v10.n2.a6},
}

\end{document}